\documentclass{amsart}

\usepackage{latexsym,amsfonts,amsmath,amssymb,amsthm,url}
\usepackage[english]{babel}
\usepackage[latin1]{inputenc}
\usepackage[T1]{fontenc}
\usepackage[dvips]{graphicx}
\usepackage[dvips]{color}

\newcommand{\R}{\mathbb{R}}

\newcommand{\N}{\mathbb{N}}
\newcommand{\E}{\mathbb{E}}
\newcommand{\Z}{\mathbb{Z}}

\newcommand{\pp}{\mathbb{P}}

\newcommand{\kF}{\mathcal{F}}

\newcommand{\kH}{\mathcal{H}}

\newcommand{\lin}{\left[\kern-0.15em\left[}
\newcommand{\rin} {\right]\kern-0.15em\right]}
\newcommand{\ilin}{\left]\kern-0.15em\left]}
\newcommand{\irin} {\right[\kern-0.15em\right[}

\newtheorem {theo} {Theorem} [section]

\newtheorem {lem} [theo] {Lemma}
\newtheorem {prop} [theo] {Proposition}
\newtheorem {cor} [theo] {Corollary}

\newtheorem {rem} [theo] {Remark}

\title[Localization on $4$ sites for VRRW.]
      {Localization on $4$ sites for Vertex-reinforced random walks on $\Z$.}

\author{Anne-Laure BASDEVANT}
\address{Laboratoire Modal'X, Universit\'e Paris Ouest Nanterre La D\'efense, B\^atiment G, 200 avenue de la R\'epublique, 92000 Nanterre, France.}
\email{anne.laure.basdevant@normalesup.org}

\author{Bruno SCHAPIRA}
\author{Arvind SINGH}
\address{D\'epartement de Math\'ematiques, B\^at. 425, Universit\'e Paris-Sud 11, F-91405 Orsay, cedex, France. }
\email{bruno.schapira@math.u-psud.fr}
\email{arvind.singh@math.u-psud.fr}


\begin{document}

\begin{abstract}
We characterize non-decreasing weight functions for which the
associated one-dimensional vertex reinforced random walk (VRRW)
localizes on $4$ sites. A phase transition appears for weights of
order $n\log \log n$: for weights growing faster than this rate, the
VRRW localizes almost surely on at most $4$ sites whereas for
weights growing slower, the VRRW cannot localize on less than $5$
sites. When $w$ is of order $n\log \log n$, the VRRW localizes
almost surely on either $4$ or $5$ sites, both events happening with
positive probability.
\end{abstract}

\maketitle

\section{Introduction}
The model of the \emph{vertex reinforced random walk} (VRRW) was
first introduced by Pemantle \cite{P} in 1992. It describes a
discrete random walk $X = (X_n,\; n\geq 0)$ on a graph $G$, which
jumps, at each unit of time $n$, from its actual position towards a
neighboring site $y$ with probability proportional to $w(Z_n(y))$,
where $w : \N \to \R_+^*$ is some deterministic weight sequence and
where $Z_n(y)$ is the local time of the walk at site $y$ and time
$n$. Thus, when $w$ is non-decreasing, the walk tends to favor sites
it has already visited many times in the past.

\medskip

A striking feature of the model is that, depending on the
reinforcement scheme $w$, it is possible for the walk to get
``trapped'' and visits only finitely many sites, even on an infinite
graph. In this case, we say that the walk \emph{localizes}. This
unusual behaviour was first observed by Pemantle and Volkov
\cite{PV} who proved that, with positive probability, the VRRW on
the integer lattice $\Z$ with linear weight $w(n) = cn+1$ visits
only $5$ sites infinitely often. This result was later completed by
Tarr\`es \cite{T1} (see also \cite{T2} for a more recent and concise
proof) who showed that localization of the walk on $5$ sites occurs
almost surely. More generally, Volkov \cite{V1} and more recently
Bena\"im and Tarr\`es \cite{BT} proved that the linearly reinforced
VRRW localizes with positive probability on any graph with bounded
degree. It is conjectured that this localization happens, in fact,
with probability $1$. However, this seems a very challenging
question as it is usually difficult to prove almost sure asymptotics
for VRRW (let us note that, even in the one-dimensional case,
Tarr\`es's proof of almost sure localization is quite elaborate).

\medskip

A seemingly closely related model is the so-called \emph{edge
reinforced random walk} (ERRW),  introduced by Coppersmith and
Diaconis \cite{CD} in 1987. The difference between VRRW and ERRW is
only that the transition probabilities for ERRW depend on the
\emph{edge} local time of the walk instead of the \emph{site} local
time. In the one-dimensional case, Davis \cite{Dav} proved that the
ERRW with non-decreasing reinforcement weight function $w$ is
recurrent (i.e. the walk visits all sites infinitely often almost
surely) i.f.f.
\begin{equation}\label{twosite}
\sum_{i=0}^{\infty}\frac{1}{w(i)} = \infty.
\end{equation}
Otherwise, the walk ultimately localizes on two consecutive sites
almost surely.

\medskip

It may seem natural to expect a similar simple criterion for VRRW.
However, the picture turns out to be much more complicated than for
ERRW because the walk may localize on subgraphs of cardinality
larger than $2$. Only partial results are currently available. For
instance, when condition
\eqref{twosite} fails and the sequence $(w(n),n\ge 0)$ is non-decreasing, the VRRW also gets stuck on two consecutive
sites\footnote{This result first appears at the end of \cite{V2}. However, there is a mistake in the original argument. For the sake of completness, we give an other proof of this result (for non-decreasing weight sequences) in Proposition \ref{proptwosite}.}. However, when \eqref{twosite} holds, the walk may or may not
localize depending on the weight function $w$. In particular, it is
conjectured that for reinforcements $w(n) \sim n^\alpha$, the walk
is recurrent for $\alpha<1$ and localizes on $5$ sites for
$\alpha=1$. In this direction, it is proved in \cite{V2} that when
$w(n)$ is of order $n^\alpha$ with $\alpha<1$,  the VRRW cannot
localize. When $\alpha<1/2$, this result was slightly refined by the
second author in \cite{Sch}, who proved that the process is either
a.s. recurrent or a.s. transient. Yet, a proof of the recurrence of
the walk in this seemingly simple setting is still missing.

\medskip

On the other hand (apart from the linear case) not much is known
about the cardinality of the set of sites visited infinitely often
when localization occurs and \eqref{twosite} holds. The aim of this
paper is to partially answer this question by investigating under
which conditions the VRRW ultimately localizes on less than $5$
sites. In order to do so, we shall associate to each weight function
$w$ a number $\alpha_c(w) \in [0,\infty]$ (the precise definition of
$\alpha_c(w)$ is given in the next section). The main result of the
paper states that:

\begin{theo}
\label{4probapositive} Assume that $w$ is non-decreasing and that
\eqref{twosite} holds. Denote by $R'$ the set of sites which are
visited infinitely often by the VRRW and by $|R'|$ its cardinality.
Then, defining $\alpha_c(w)$ as in \eqref{defAlphac}, we have
\begin{eqnarray*}
|R'|=4 \textrm{ with positive probability} &\Longleftrightarrow & \alpha_c(w)<\infty.\\
|R'|=4 \textrm{ almost surely} &\Longleftrightarrow & \alpha_c(w)=0.\\
\begin{array}{c}
|R'| \textrm{ equals $4$ or $5$ a.s., both events }\\
\textrm{occurring with positive probability}
\end{array}  &\Longleftrightarrow & \alpha_c(w)\in (0,\infty).
\end{eqnarray*}
\end{theo}
It is easy to see that a VRRW can never localize on $3$ sites (see
Proposition \ref{R'3}) therefore Theorem \ref{4probapositive}
combined with criterion \eqref{twosite} for localization on two sites covers all
possible cases where a VRRW with non-decreasing weight localizes with
positive probability on less than  $5$ sites.

\medskip

The last part of the theorem shows that the size of $R'$ can itself
be random. Such a result was already observed for graphs like
$\Z^d$, $d\ge 2$ (for linear reinforcement) where different
non-trivial localization patterns may occur (see \cite{BT,V1}). Yet
we find this result more surprising in the one-dimensional setting
since $R'$ is necessarily an interval.

\medskip

The parameter $\alpha_c(w)$ can be explicitly calculated for a large
class of weights $w$. In particular, if $w(n)\sim n\log\log n$, then $\alpha_c(w) = 1$. Moreover,
\begin{prop}
\label{alphac_cor} For any non-decreasing weight function $w$ such
that \eqref{twosite} holds:
\begin{eqnarray*}
w(n)\ \gg  \footnotemark \ n\log\log n & \Longrightarrow & \alpha_c(w) =
0. \\
w(n)\ \asymp \footnotemark\
n\log\log n & \Longrightarrow & \alpha_c(w) \in (0,\infty). \\
w(n)\ \ll \ n\log\log n & \Longrightarrow & \alpha_c(w) = \infty.
\end{eqnarray*}
\addtocounter{footnote}{-1}
\footnotetext{we use the notation $f\gg g$, when $f(n)/g(n)\to \infty$.}
\stepcounter{footnote}
\footnotetext{we say that $f\asymp g$ when there exists a constant $c>0$, such that
$c^{-1}\, f(n) \le g(n) \le c\, f(n)$, for all $n$ large enough.}
\end{prop}

Let us mention that, when $\alpha_c(w) = \infty$, Theorem
\ref{4probapositive} simply states that, if the walk localizes,
$|R'| \geq 5$ necessarily. In fact, it is proved in a forthcoming
paper \cite{BSS} that there exist non-decreasing weight functions
$w$ for which the VRRW localizes almost surely on finite sets of
arbitrarily large cardinality (this result is, in a way, similar to
those proved in \cite{ETW1,ETW2} for another related model of
self-interacting random walks).

\medskip

The proof of Theorem \ref{4probapositive} is based on two main
techniques. First we use martingales arguments which were introduced
by Tarr\`es in \cite{T1,T2}. These martingales have the advantage of
taking into account the facts that, on each site, the process is
roughly governed by an urn process, but also the fact that all these
urns are strongly correlated. The second tool is a continuous time
construction of the VRRW, called Rubin's construction, which was
already used by Davis \cite{Dav} for urn processes, and by Sellke in
\cite{Sel} in the case of edge reinforcement. Tarr\`es introduced in
\cite{T2} a variant of this construction, which allows for powerful
couplings in the case of non-decreasing weights and which will be
very useful in this study.

\medskip
The remainder of the paper is organized as follows. In the next
section we give some simple results concerning $w$-urns processes
which will play an important role in the proof of the theorem. In
Section \ref{section3}, we recall some classical results concerning
VRRW. The proof of Theorem \ref{4probapositive} is provided in
Section \ref{section4}. Finally we prove Proposition
\ref{alphac_cor} in the appendix along with other technical lemmas
concerning properties of the critical parameter $\alpha_c(w)$.

\section{$w$-urn processes}
\label{section2}

\subsection{Weight function $w$ and the parameter $\alpha_c$}
In the rest of the paper, we call \emph{weight sequence} $w$ a sequence $(w(n),n\ge 0)$ of positive real numbers. It will be convenient to extend $w$ into a weight
function $w: \R_+ \to \R_+^*$ by $w(t) = w(\lfloor t\rfloor)$ where $\lfloor t\rfloor$
stands for the integer part of $t$. Then, given $w$, we set
\begin{equation*}
W(t) := \int_0^t \frac 1{w(u)}\, du.
\end{equation*}
When condition (\ref{twosite}) holds, we have $W(\infty) = \infty$
and $W$ is an homeomorphism of $\R_+$ whose inverse we denote by
$W^{-1}$. Then, for $\alpha>0$, we define the integral
\begin{equation}\label{defI}
I_\alpha(w):= \int_0^\infty
\frac{dx}{w(W^{-1}(W(x)+\alpha))}=\int_0^\infty
\frac{w(W^{-1}(y))}{w(W^{-1}(y+\alpha))}\, dy.
\end{equation}
If furthermore we assume that $w$ is non-decreasing, then $\alpha
\to I_\alpha(w)$ is non-increasing and we can define the critical
parameter $\alpha_c(w)$ by
\begin{equation}\label{defAlphac}
\alpha_c(w):=\inf \{\alpha\ge 0 \ :\ I_\alpha(w)<\infty\}\in
[0,\infty]
\end{equation}
with the convention that $\inf \emptyset = \infty$.

\subsection{$w$-urn processes}
A $w$-urn is a process $(R_n,B_n)_{n\ge 0}$ defined on some probability space $(\Omega,\mathcal{F},\pp)$, such that for all $n\ge
0$, $R_n+B_n=n$, $R_{n+1}\in \{R_n,R_n+1\}$, and
\begin{equation*}
\pp\{R_{n+1}=R_n + 1\} = \frac {w(R_n)}{w(R_n)+w(B_n)}.
\end{equation*}
 We call $R_n$ (resp. $B_n$) the number of red (resp. blue)
balls in the urn after the $n$-th draw. Set $R_\infty =
\lim_{n\to \infty}R_n$ and $B_\infty=\lim_{n\to \infty} B_n$.

Our interest towards $w$-urn processes comes from fact that, if we
consider a VRRW on the finite set $\{-1,0,1\}$ (\emph{i.e.} the walk
reflected at $1$ and $-1$, see Section \ref{section3.4}), then joint
local times of the walk at sites $1$ and $-1$ and at time $2n$ is
exactly a $w$-urn process. The next proposition describes the
asymptotic behavior of such an urn. Several arguments used during
the proof of the result below will also play an important role when
proving Theorem \ref{4probapositive}.

\begin{prop}
\label{propurne} For any weight sequence $w$ (not necessarily
non-decreasing), we have
\begin{eqnarray*}
\sum_{n\ge 0} \frac 1 {w(n)}< +\infty \quad \Longleftrightarrow
\quad R_\infty<+\infty\ \textrm{or}\ B_\infty <+\infty\;\hbox{ a.s.}
\end{eqnarray*}
The process $\hat{M} = (\hat{M}_n,\, n\geq 0)$ defined by $\hat{M}_n
:=W(R_n)-W(B_n)$ is a martingale. Moreover,
\begin{enumerate}
\item[(a)] If $\sum_n 1/w(n)=\infty$ and $\sum_n 1/w(n)^2<+\infty$, then $\hat{M}_n$ converges a.s. to some random variable $\hat{M}_\infty$, which
admits a symmetric density with unbounded support.
\item[(b)] If $\sum_n 1/w(n)^2 = \infty$ and $\inf_n w(n) >0$, then
$\liminf \hat{M}_n = -\infty$ and $\limsup \hat{M}_n=+\infty$ a.s.
\end{enumerate}
\end{prop}
\begin{proof}
The first equivalence is well known (see, for instance \cite{P2})
and follows immediately from Rubin's construction of the urn process
which we will recall below. Let us note that we can rewrite $\hat{M}$ in the form
$$
\hat{M}_n = \sum_{k=0}^{n-1} \left(\frac {1_{\{\textrm{the
$(k+1)$-th draw is Red}\}}}{w(R_k)} - \frac {1_{\{\textrm{the
$(k+1)$-th draw is Blue}\}}}{w(B_k)}\right).
$$
Therefore, $\hat{M}$ is clearly a martingale. Let
\begin{eqnarray*} V_n &:=& \sum_{k\le n} (\hat{M}_{k+1}-\hat{M}_k)^2\\
  &=& \sum_{k\le n } \left(\frac {1_{\{\textrm{the $(k+1)$-th draw is Red}\}}}{w(R_k)^2} + \frac {1_{\{\textrm{the $(k+1)$-th draw is
  Blue}\}}}{w(B_k)^2}\right).
\end{eqnarray*}
We have
\begin{equation*}
\E[\hat{M}_n^2]  = \E[ V_{n-1} ] \leq 2\sum_{k =
0}^{\infty}\frac{1}{w(k)^2}.
\end{equation*}
Thus, when assumption (a) holds, $\hat{M}$ converges almost surely
and in $L^2$ towards some random variable $\hat{M}_\infty$. We now
use Rubin's construction to identify this limit: let $(\xi_n,n\ge
0)$ and $(\xi'_n,n\ge 0)$ be two sequences of independent
exponential random variables with mean $1$. Define the random times
$t_k=\xi_0/w(0) + \dots +\xi_k/w(k)$ and $t'_k=\xi'_0/w(0)+\dots
+\xi'_k/w(k)$. We can construct the $w$-urn process $(R_n,B_n)$ from
these two sequences by adding a red ball in the urn at each instant
$(t_k)_{k\geq 0}$ and a blue ball at each instant $(t'_k)_{k\geq 0}$
(see the appendix in \cite{Dav} for details). Using this
construction, we can rewrite $\hat{M}_n$ in the form
$$\hat{M}_n=\sum_{k=0}^{R_n-1} \frac {1-\xi_k} {w(k)}- \sum_{k=0}^{B_n-1}  \frac {1-\xi_k'} {w(k)}+
\sum_{k=0}^{R_n-1}  \frac {\xi_k} {w(k)}-\sum_{k=0}^{B_n-1}  \frac
{\xi'_k} {w(k)}.$$ Observe that, by construction, for any $n$,
$$\left|\sum_{k=0}^{R_n-1}  \frac {\xi_k} {w(k)}-\sum_{k=0}^{B_n-1}  \frac {\xi'_k} {w(k)}\right|
 \le \max\left( \frac {\xi_{R_{n}}} {w(R_n)},  \frac {\xi'_{B_{n}}} {w(B_n)}\right).$$
Since a.s. $R_n\wedge B_n \to \infty$, we deduce that the r.h.s.
above tends to $0$ a.s. hence
\begin{equation}\label{limitUrn}
\hat{M}_\infty = \sum_{k=0}^{\infty} \frac {1-\xi_k} {w(k)}-
\sum_{k=0}^{\infty}  \frac {1-\xi_k'} {w(k)}.
\end{equation}
 Both sums in the r.h.s. of the previous
equation converge because they have a finite second moment. Thus,
$\hat{M}_\infty$ admits a symmetric density with unbounded support
since the $\xi_n$'s and $\xi'_n$'s are independent and $\xi_0$ has a
non-vanishing density on $\R_+$.

It remains to prove $(b)$. Let us observe that, when $\inf_n
w(n)>0$, the martingale $\hat{M}_n$ has bounded increments. Using
Theorem 2.14 in \cite{HH}, it follows that, a.s., either $\hat{M}_n$
converges or $\limsup \hat{M}_n = +\infty$ and $\liminf \hat{M}_n =
-\infty$. Moreover, when $\sum_n 1/w(n)^2 = \infty$, we have
$\lim_{n\to\infty} V_n = \infty$ a.s. Therefore, we can define
$k_n:=\inf\{m\ :\ V_m\ge n\}$ and Theorem $3.2$ of \cite{HH} states
that $\hat{M}_{k_n}/\sqrt n$ converges in law towards a standard
normal variable. In particular, $\hat{M}$ cannot converge. This
completes the proof of the proposition.
\end{proof}

The next result illustrates how the parameter $\alpha_c(w)$ of
Theorem \ref{4probapositive} naturally appears in connection with
$w$-urns.
\begin{cor}\label{Yurn} Consider a $w$-urn $(R_n,B_n)$. Assume that $w$ is non-decreasing with $\sum_n
1/w(n)=\infty$ and $\sum_n 1/w(n)^2 <\infty$ and set
\begin{equation*}
Y^B := \sum_{k=0}^\infty \frac{1_{\{\textrm{the $(k+1)$-th draw is
Blue}\}}}{w(k)},\qquad Y^R := \sum_{k=0}^\infty \frac{1_{\{\textrm{the
$(k+1)$-th draw is Red}\}}}{w(k)}.
\end{equation*}
Then, we have
\begin{itemize}
\item[(i)]If $\alpha_c(w)=0$ then, a.s.,  $\min(Y^B,Y^R)<\infty$.
\item[(ii)]If $\alpha_c(w)\in (0,\infty)$ then $\pp\{\min(Y^B,Y^R)<\infty\}\in (0,1)$.
\item[(iii)]If $\alpha_c(w)=\infty$  then, a.s.,
$\min(Y^B,Y^R)=\infty$.
\end{itemize}
\end{cor}

\begin{proof}
 According to Proposition \ref{propurne}, $W(R_n)-W(B_n)$ converges to some
 random variable $\hat{M}_\infty$  with a symmetric density and unbounded
support. Let $\delta=|\hat{M}_\infty|/2$. On the event
$\{\hat{M}_\infty>0\}$, we have, for $n$ large enough,
$$W(n)\ge W(R_n) \ge W(B_n)+ \delta$$
which yields, for some (random but finite) constant $c$
$$Y^B\le c\sum_{k=0}^\infty \frac{1_{\{\textrm{the $(k+1)$-th draw is Blue}\}}}{w(W^{-1}(W(B_k)+ \delta))}= \sum_{k=0}^\infty \frac{c}{w(W^{-1}(W(k)+ \delta))}.$$
Thus, by symmetry and using $\pp\{\hat{M}_\infty =0\} = 0$, we get,
a.s.,
$$\min(Y^B,Y^R) \leq \sum_{k=0}^\infty \frac{c}{w(W^{-1}(W(k)+ \delta))}.$$
This proves (i) and also that $\pp\{\min(Y^B,Y^R)<\infty\}>0$
whenever $\alpha_c(w)\in (0,\infty)$. Conversely, set $\delta' =
2|\hat{M}_\infty|$. On the event $\{\hat{M}_\infty \geq 0\}$, we
have, for $n$ large enough
$$n=B_n+R_n \le B_n+ W^{-1}(W(B_n)+ \delta') \le R_n+ W^{-1}(W(R_n)+ \delta').$$
This gives
$$Y^B\ge c\sum_{k=0}^\infty \frac{1_{\{\textrm{the $(k+1)$-th draw is Blue}\}}}{w(B_k+W^{-1}(W(B_k)+ \delta'))}= \sum_{k=0}^\infty \frac{c}{w(k+W^{-1}(W(k)+ \delta'))}.$$
and the same bound also holds for $Y^R$. Therefore, by symmetry, we
get, a.s.,
$$\min(Y^B,Y^R) \geq \sum_{k=0}^\infty \frac{c}{w(k+W^{-1}(W(k)+ \delta'))}.$$
We conclude the proof using Lemma \ref{teclem} of the appendix which
insures that the sum above is infinite when $\delta' <
\alpha_c(w)$.
\end{proof}

\section{Vertex reinforced random walk}
\label{section3}

\subsection{The VRRW} In the remainder of the paper, given the weight sequence $w$, $(X_n,n\ge 0)$
will denote a nearest neighbour random walk on the integer lattice
$\Z$, starting from $X_0=0$ with transition probabilities given by
\begin{equation}\label{transiVRRW}
\pp\{X_{n+1}= x\pm 1 \mid \kF_n\}\, = \, \frac{w(Z_n(x\pm
1))}{w(Z_n(x+1)) + w(Z_n(x-1))},
\end{equation}
 where $(\kF_n,n\ge 0)$ is the natural filtration $\sigma(X_0,\ldots,X_n)$ and
$Z_n(y)$ stands, up to a constant, for the local time of $X$ at site $y$ and at time
$n$:
$$Z_n(y):=z_0(y)+\sum_{k=0}^n 1_{\{X_k=y\}}.$$
We call the sequence $\mathcal{C} := (z_0(y),\, y\in\Z)$ the\emph{
initial local time configuration}. We say that $X$ is a VRRW when it
starts from  the trivial configuration $\mathcal{C}_0 :=
(0,0,\ldots)$. However, it will sometimes be convenient to consider
the walk starting from some other configuration $\mathcal{C}$. In
that case, we shall mention it explicitly and emphasize  this fact
by calling $X$ a $\mathcal{C}$-VRRW.

\medskip

In the rest of this section, we collect some important results
concerning the VRRW which we will use during the proof of Theorem
\ref{4probapositive} in Section \ref{section4}. For additional
details, we refer the reader to \cite{Dav,P2,Sel,T1,T2} and the
references therein.

\subsection{The martingales $M_n(x)$}\label{section3.2}

For $x\in \Z$, define $Z_\infty(x) := \lim_{n\to \infty} Z_n(x)$.
Recall that $R'$ stands for the set of sites visited infinitely
often by the walk:
$$R' :=\{x\in \Z\ :\ Z_\infty(x)=\infty\}.$$
The following quantities will be of interest:
\begin{equation}\label{defY}
Y_n^{\pm}(x) := \sum_{k=0}^{n-1} \frac{1_{\{X_k=x\textrm{ and
}X_{k+1}=x\pm 1\}}}{w(Z_k(x\pm 1))},
\end{equation}
 and
$$M_n(x) := Y_n^+(x)-Y_n^-(x).$$
It is a basic observation due to Tarr\`es \cite{T1,T2} that
$(M_n(x),n\ge 1)$ is a martingale for each $x\in \Z$. Moreover, if
\begin{equation}\label{wsquare}
\sum_{n=0}^{\infty} \frac{1}{w(n)^2}<\infty,
\end{equation}
then these martingales are bounded in $L^2$, and thus converge a.s.
and in $L^2$ towards
$$
M_\infty(x) := \lim_{n\to\infty} M_n(x).
$$
We will also consider the (possibly infinite) limits:
$$
Y_\infty^\pm(x):= \lim_{n\to \infty} Y_n^{\pm}(x).
$$
From the definition of $Y^{\pm}$, we directly obtain the identity
\begin{eqnarray}
 \label{Ynpm}
Y_n^+(x-1) + Y_n^-(x+1) = W(Z_{n}(x))-W(1)1_{\{x=0\}},
\end{eqnarray}
which holds for all $x\in \Z$ and all $n\ge 0$. In particular, we
get
\begin{multline*}
W(Z_{n}(x+2))-W(Z_{n}(x)) = Y_n^-(x+3)-Y_n^+(x-1)+M_{n}(x+1)\\
+W(1)(1_{\{x=-2\}}-1_{\{x=0\}}).
\end{multline*}
More generally, if we now consider a $\mathcal{C}$-VRRW starting
from some arbitrary initial local time configuration $\mathcal{C}$,
then $M_n(x)$ is still a martingale and the equation above takes the
form
\begin{equation}\label{eqW}
W(Z_{n}(x+2))-W(Z_{n}(x))=Y_n^-(x+3)-Y_n^+(x-1)+M_{n}(x+1)+
c(x,\mathcal{C}).
\end{equation}
where $c(x,\mathcal{C})$ is some constant depending only on $x$ and
on the configuration $\mathcal{C}$.

\subsection{Time-line construction of the VRRW}
We now describe a method to construct the VRRW for a collection of
exponential random variables which is in a way similar to Rubin's
algorithm for $w$-urns. This construction was introduced by Tarrès
in \cite{T2} and may be seen as a variant for the VRRW of the
continuous time construction previously described by Sellke in
\cite{Sel} for edge reinforced random walks. One of the main
advantages of this construction is that it enables to create
non-trivial coupling between VRRWs.
\begin{figure}
\begin{center}
\includegraphics[width=12.6cm,angle=0]{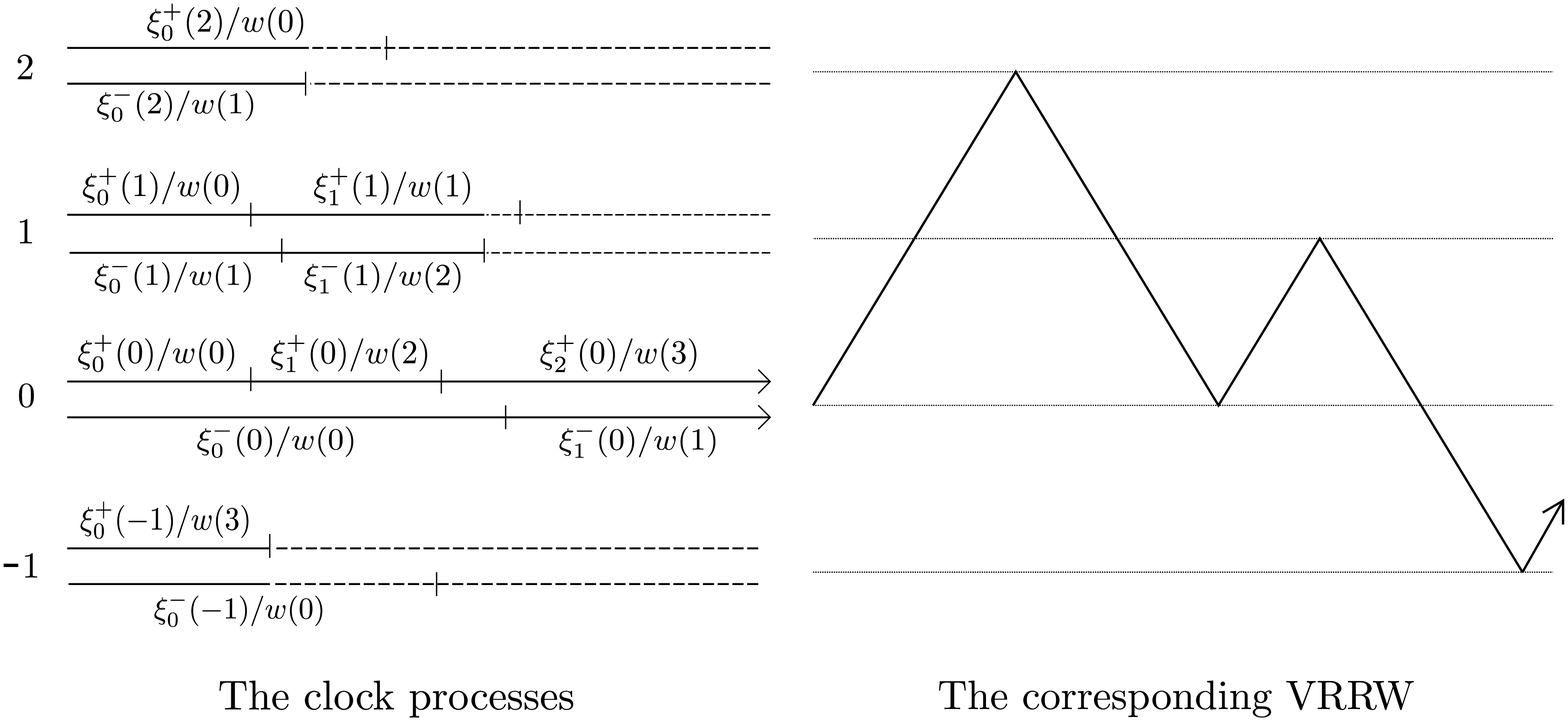}
\end{center}
\caption{Illustration of the time-line construction.}
\end{figure}
Let us fix a sequence $$\xi := (\xi_n^\pm(y),n\ge 0, y\in\Z) \in \R_+^\N$$ of positive
real numbers. The value $\xi^-_n(y)$ (resp. $\xi^+_n(y)$) will be
related to the duration of a clock attached to the oriented edge
$(y,y-1)$ (resp. $(y,y+1)$). Given this sequence, we create a
deterministic, integer valued, continuous-time process $(\widetilde
X(t),t\ge 0)$ in the following way:
\begin{itemize}
\item Set $\widetilde X(0)=0$ and  attach two clocks to the oriented edges  $(0,-1)$ and
$(0,1)$ ringing respectively at times $\xi_0^-(0)/w(0)$ and
$\xi_0^+(0)/w(0)$.
\item When the first clock rings at time $\tau_1:=\xi_0^+(0)/w(0)\wedge \xi_0^-(0)/w(0)$, stop both clocks
and set $\widetilde X(\tau_1)=\pm1$ depending on which clock rung
first. (if both clocks ring at the same time, we decide that
$\widetilde X$ stays at $0$ forever).
 \end{itemize}
Assume that we have constructed $\widetilde{X}$ up to some time $t>0$ at
which time the process makes a right jump from some site $x-1$ to $x$. Denote by $k$
the number of jumps from $x$ to $x-1$ and  by $m$ the number of visits
to $x-1$ before time $t$. We follow the procedure below:
 \begin{itemize}
\item Start a new clock attached to the oriented edge $(x,x-1)$, which will ring
after a time $\xi_k^-(x)/w(m)$.
\item If the process already visited $x$ some time in the past, restart the clock attached to the oriented edge $(x,x+1)$ which had previously
been stopped when the process last left site $x$. Otherwise, start the first clock for this edge which will ring at time $\xi_0^+(x)/w(0)$.
\item As soon as one of these two clocks rings, stop both of them and let the process jump along the edge corresponding to the clock which rung
first (if both clocks ring at the same time, we decide
that $\widetilde X$ stays in $x$ forever).
 \end{itemize}
We use a similar rule when the process makes a left jump from some
site $x$ to $x-1$. We say that this construction \emph{fails} if at
some time, two clocks ring simultaneously. Let now $\tau_i$ stand
for the time of the $i$-th jump of $\widetilde{X}$ (with the
convention $\tau_0=0$ and $\tau_{n+1}=\tau_{n}$ if $\widetilde X$
does not move after time $\tau_n$) and define the discrete time
process $X = (X_n,\, n\geq 0)$ by
\begin{equation*}
X_n := \widetilde{X}(\tau_n).
\end{equation*}
It is an elementary observation that if we now choose the
$\xi_n^\pm(y)$ to be independent exponential random variables with
mean $1$, then the construction does not fail with probability $1$
and the resulting process $X$ is a VRRW with weight $w$.

\begin{rem} For the sake of clarity,  we only describe the
construction for the VRRW starting from the trivial configuration $\mathcal{C}_0$.
However, it is clear that we can do a similar construction for any $\mathcal{C}$-VRRW  by simply
replacing the duration of the clocks $\xi_k^\pm(x)/w(m)$ with $\xi_k^\pm(x)/w(z_0(x\pm 1) + m)$.
\end{rem}

 A remarkable feature of this construction comes from the fact
that we can simultaneously create a family $(\widetilde{X}^{(u)},\;
u\geq 0)$ of processes with nice monotonicity properties with
respect to the $u$ parameter. To this end, define, for $x\in\Z$,
$$\kH_x:= \left((\xi_n^\pm(y),n\ge 0, y\neq x),(\xi_n^-(x),n\ge 0),(\xi_n^+(x),n\ge 1)\right) \in \R_+^\N.$$
Then, given $\kH_x$ together with a real number $u>0$, the pair
$(\kH_x,u)$ defines a deterministic process $X^{(u)} =
(X_n^{(u)},n\geq 0)$ using the construction above with $\xi_0^+(x) =
u$. The following lemma is easily obtained by induction.
\begin{lem}[Tarr\`es \cite{T2}]
\label{tarres1} Suppose that $w$ is non-decreasing. Fix $\kH_x$ and
$0 < u\le u'$ and assume that the construction for $X^{(u)}$ and
$X^{(u')}$ both succeed. Given $y\in \Z$ and $k\ge 1$, denote by
$\sigma$ (resp. $\sigma'$) the time when $X^{(u)}$ (resp.
$X^{(u')}$) visits $y$ for the $k$-th time. If $\sigma$ and
$\sigma'$ are both finite, then
\begin{eqnarray*}
Z_\sigma^{(u)}(y+1) \ge Z_{\sigma'}^{(u')}(y+1) \quad &\textrm{and}& \quad Z_\sigma^{(u)}(y-1) \le Z_{\sigma'}^{(u')}(y-1)\\
N_\sigma^{(u)}(y,y+1) \ge N_{\sigma'}^{(u')}(y,y+1) \quad
&\textrm{and}& \quad N_\sigma^{(u)}(y,y-1) \le
N_{\sigma'}^{(u')}(y,y-1),
\end{eqnarray*}
where $Z_n^{(s)}$ stands the local time of $X^{(s)}$ and
$N_n^{(s)}(y,y\pm 1)$ denotes the number of jumps from $y$ to $y\pm
1$ up to time $n$. Moreover, denote by $\theta^{\pm}$ (resp.
${\theta'}^{\pm}$) the time when  $X^{(u)}$ (resp. $X^{(u')}$) jumps
for the $k$-th time from $y$ to $y\pm1$. If these quantities are
finite, then
\begin{equation*}
Y^{(u)+}_{\theta^{+}}(y) \le
Y^{(u')+}_{{\theta'}^{+}}(y)\qquad\hbox{and}\qquad
Y^{(u)-}_{\theta^{-}}(y) \ge Y^{(u')-}_{{\theta'}^{-}}(y),
\end{equation*}
where $Y^{(s)\pm}$ is defined as in \eqref{defY} for the process $X^{(s)}$.
\end{lem}

The combination of the time-line construction of the walk from
i.i.d. exponential random variables together with Lemma
\ref{tarres1} yields a simple proof of the following key result
concerning the localization of the VRRW:
\begin{lem}[Tarr\`es \cite{T2}]
\label{Yfini} Assume that $w$ is non-decreasing and that $\sum_n
1/w(n)^2$ is finite. Then, for any $x\in \Z$, a.s.,
$$\{Y_\infty^+(x)<\infty\}=\{Y_\infty^-(x)<\infty\}\ = \ \{Z_\infty(x-1)<\infty\} \cup \{Z_\infty(x+1)<\infty\}.$$
\end{lem}
\begin{proof} This result is proved in \cite{T2} only for linear reinforcements $w$ but
the same arguments apply, in fact, for any non-decreasing weight
function. However, since some details are omitted in \cite{T2}, for
the sake of completeness, we provide here a detailed proof
(differing in some aspects from the original one). The first
equality $\{Y_\infty^+(x)<\infty\}=\{Y_\infty^-(x)<\infty\}$ follows
from the fact that the martingale $M_n(x)$ converges a.s. to some
finite limit when $\sum_n 1/w(n)^2$ is finite. Concerning the second
equality, the inclusion
$$\{Y_\infty^+(x)<\infty\}=\{Y_\infty^-(x)<\infty\}\ \supset \ \{Z_\infty(x-1)<\infty\} \cup \{Z_\infty(x+1)<\infty\}.$$
is straightforward (one of the sums $Y_\infty^{\pm}$ has only a
finite number of terms). We use the time-line construction of the
VRRW $X$ from the sequence $(\xi_n^\pm(y),n\ge 0,y\in \Z)$ to prove
the converse inclusion. Denote by $N_k(x, x\pm 1)$ the number of
jumps of $X$ from $x$ to $x\pm 1$ before time $k$, and set
$$T_x^\pm :=\sum_{k\ge 0} \frac{1_{\{X_k=x,\ X_{k+1}=x\pm 1\}}\xi_{N_k(x, x\pm 1)}^\pm(x)}{w(Z_k(x\pm 1))}.$$
Thus, $T_x^\pm$ represents the total time consumed by the
clocks attached to oriented edge $(x,x\pm 1)$. We claim that
\begin{equation}
\label{Tx} \{Y_\infty^+(x)<\infty\}\cap \{Z_\infty(x-1)=\infty\}\cap
\{Z_\infty(x+1) =\infty\} \subset \{T_x^+=T_x^-<\infty\}.
\end{equation}
We prove the result for $x<0$ (the proof for $x>0$ and $x=0$ are
similar). Let $\theta_k$ denote the time of the $k$-th jump from
site $x+1$ to $x$ and let $i_k$ be the local time at site $x+1$ and
at time $\theta_k$. With this notation, on the event
$\{Z_\infty(x+1) =\infty\}$, we can write
\begin{equation*}
T_x^+= \sum_{k\ge 1} \xi_{k-1}^+(x)\frac{1_{\{\theta_k <
\infty\}}}{w(i_k)}\quad \mbox{ and } \quad Y_\infty^+(x)= \sum_{k\ge
1} \frac{1_{\{\theta_k < \infty\}}}{w(i_k)}.
\end{equation*}
Define now
\begin{equation*}
T_n^+(x):= \sum_{k=1}^{n-1} (\xi_{k-1}^+(x)-1)\frac{1_{\{\theta_k <
\infty\}}}{w(i_k)}.
\end{equation*}
Recall that $(\kF_n,n\ge 0)$ stands for the natural filtration of
$X$ and notice that $i_k$ is $\kF_{\theta_k}$-measurable whereas
$\xi_{k-1}^+(x)$ is independent of $\kF_{\theta_k}$. Thus,  $(T_n^+(x),n\ge 1)$ is a $\kF_{\theta_n}$-martingale. Moreover, using  that $w(i_k)\ge w(k)$ for all $k\ge 1$, it follows
that the $L^2$-norm of this martingale is bounded by $\sum_{k\ge 0} {w(k)^{-2}}<\infty$.
In particular, this implies that
$Y_\infty^+(x)$ is finite if and only if $T_x^+$ is finite. It is
also clear from the construction of the time-line process that, on
the event $\{Z_\infty(x-1)=\infty\}\cap \{Z_\infty(x+1) =\infty\}$,
we have $T_x^+=T_x^-$. Thus we have established \eqref{Tx}.
\medskip

It remains to prove that the event $$\mathcal{E}_x:=\{T_x^+=T_x^- <\infty\}\cap \{Z_\infty(x-1)=\infty\}\cap
\{Z_\infty(x+1) =\infty\}$$ has
probability $0$.
Recall the notation
$$\kH_x:= \left((\xi_n^\pm(y),n\ge 0)_{y\neq x},(\xi_n^-(x),n\ge 0),(\xi_n^+(x),n\ge 1)\right) \in \R_+^\N.$$
and denote by $\mu$ the product measure on $\R_+^\N$ under which
$\kH_x$ is a collection of i.i.d. exponential random variables with
mean $1$. Given $\xi_0^+(x)$, the pair $(\kH_x,\xi_0^+(x))$ defines
the  (deterministic) process $X = X(\kH_x,\xi_0^+(x))$ via the
time-line construction and $X$ under the product law
$\pp:=\mu\times\hbox{Exp}(1)$ is a VRRW.

Let us note that, for $\mu$-a.e. realization of $\kH_x$, the set of values of $\xi_0^+(x)$
such that the time-line construction fails  is countable hence has zero Lebesgue measure.
 Moreover, Lemma \ref{tarres1} implies that, for any $(\kH_x,u)$ and $(\kH_x,u')$ in $\mathcal{E}_x$ with $u'>u$, we have
\begin{equation*}
T_x^+(\kH_x,u') > T_x^+(\kH_x,u)\quad\hbox{ and }\quad
T_x^-(\kH_x,u') \leq T_x^-(\kH_x,u).
\end{equation*}
Thus, for any $\kH_x$, there is at most one value of $\xi_0^+(x)$ such that $(\kH_x,\xi_0^+(x))\in \mathcal{E}_x$. This yields
$$\pp\{\mathcal{E}_x\}=\mathbb{E}_\mu\left(\int_0^\infty e^{-u} 1_{\{(\kH_x,u)\in\mathcal{E}_x\}}du\right)=\mathbb{E}_\mu(0)=0.$$
\end{proof}

A weaker statement can also be obtained when the assumptions of
Lemma \ref{Yfini} do not hold.
\begin{lem}\label{Yfini2}
For any weight sequence $w$ and for any $x\in \Z$, we have, a.s.
$$ \{Z_\infty(x-1)<\infty\} \cup \{Z_\infty(x+1)<\infty\}\subset \{Y_\infty^+(x)<\infty\}\cap\{Y_\infty^-(x)<\infty\}.$$
\end{lem}

\begin{proof}
By symmetry, we can assume without loss of generality that
$Z_\infty(x-1)<\infty$. On the one hand,  $Y_\infty^-(x)$ is finite
since it is a sum with a finite number of terms. On the other hand,
the conditional Borel-Cantelli lemma implies that
$Y_\infty^+(x)<\infty$ (apply for instance the theorem of \cite{Ch}
with the sequence $1_{\{X_k=x-1\}}$).
\end{proof}

\subsection{VRRW restricted to a finite set.}\label{section3.4}
In the sequel, it will be convenient to consider the vertex
reinforced random walk restricted to some interval $\lin a,b\rin:=
\{x\in \Z\ :\ a\le x\le b\}$ for some $a \leq 0 \leq b$,
\emph{i.e.} a walk with the same transition probabilities
\eqref{transiVRRW} as the VRRW $X$ on $\Z$ except at the boundary
sites $a$ and $b$ where it is reflected. We shall use the notation
$\bar{X}$ to denote this reflected process. We also add a bar to denote all the quantities
$\bar{Z}$,$\bar{Y}^\pm$,$\bar{M}$,\ldots related with the reflected
process $\bar{X}$.

\begin{rem} Let us emphasize the fact that, for $x\in \ilin a,b \irin$, the processes $\bar{M}_n(x):=\bar{Y}^{+}_n(x)-\bar{Y}^{-}_n(x)$ are
still martingales, which are bounded in $L^2$ when \eqref{wsquare}
holds. In particular,  Lemma \ref{Yfini} still holds for the
reflected random walk for all site $x\in \ilin a,b \irin$.
 However,  $\bar{M}_n(a)$ and
$\bar{M}_n(b)$ are not martingales anymore. In particular,
$\bar{Y}^+_\infty(a)$ or $\bar{Y}^{-}_\infty(b)$ can be infinite
whereas $\bar{Y}^-_\infty(a)$ and $\bar{Y}^{+}_\infty(b)$ are, by
construction, always equal to $0$.
\end{rem}

We can construct the VRRW $\bar{X}$
restricted to the interval $\lin a,b\rin$ using the same time-line
construction used for $X$, choosing again the random variables
$\xi_k^\pm(x)$ independent and exponentially distributed except for
the two boundary r.v. $\xi_0^-(a)$ and $\xi_0^+(b)$ which are now chosen
equal to $\infty$ (this prevent the walk from ever jumping from $a$
to $a-1$ or from $b$ to $b+1$). Let us note that, this construction
depends only upon $(\xi^\pm_n(x), n\geq 0, x\in \ilin a,b\irin )$.
\medskip

Let $\bar{X}'$ denote another VRRW restricted to $\lin a,b'\rin
\supset \lin a,b\rin$ for some $b'\geq b$. Using the time-line
construction for $\bar{X}$ and $\bar{X}'$ with the same random
variables $\xi_k^\pm(x)$, except for $\xi_0^+(b)$, we directly deduce from
Lemma \ref{tarres1} a monotonicity result between the local time processes of
$X$ and $X'$ :

\begin{lem}\label{couplrest}
 Assume that $w$ is non-decreasing. Fix $z\in \lin a,b\rin$ and $k\ge 1$,
 let $\sigma,\sigma'$ be the times when $\bar{X},\bar{X}'$ visit $z$
for the $k$-th times. On the event $\{\sigma <\infty \hbox{ and } \sigma'<\infty\}$, we have, a.s.,
\begin{eqnarray*}
\bar{Z}_\sigma(z+1) \le \bar{Z}'_{\sigma'}(z+1) \quad &\textrm{and}& \quad \bar{Z}_\sigma(z-1) \ge \bar{Z}'_{\sigma'}(z-1)\\
\bar{N}_\sigma(z,z+1) \le \bar{N}'_{\sigma'}(z,z+1) \quad
&\textrm{and}& \quad \bar{N}_\sigma(z,z-1) \ge
\bar{N}'_{\sigma'}(z,z-1),
\end{eqnarray*}
where $N$ and $N'$ are defined as in Lemma \ref{tarres1} for
$\bar{X}$ and $\bar{X}'$. Moreover, if we denote by $\theta^{\pm}$
(resp. ${\theta'}^{\pm}$) the time when $\bar{X}$ (resp. $\bar{X}'$)
jump for the $k$-th time from $z$ to $z\pm1$, then, on the event of
these quantities being finite, we have, a.s.,
\begin{equation*}
\bar{Y}^{+}_{\theta^{+}}(z) \ge
\bar{Y}'^{+}_{{\theta'}^{+}}(z)\qquad\hbox{and}\qquad
\bar{Y}^{-}_{\theta^{-}}(z) \le \bar{Y}'^{-}_{{\theta'}^{-}}(z).
\end{equation*}
\end{lem}

We conclude this section with a simple lemma we will
repeatedly invoke to reduce the study of the localization properties
of the VRRW $X$ on $\Z$ to those of the VRRW $\bar{X}$
restricted to a finite set.
\begin{lem}\label{lemrest2}
 For $N>0$, let $\bar{X}$ be a VRRW on
$\lin 0,N\rin$ and define the events
\begin{eqnarray*}
\mathcal{E}&=&\{\bar{Y}^+_\infty(0)<\infty\}\cap\{\bar{Y}^{-}_\infty(N)<\infty\},\\
\mathcal{E}' &=& \mathcal{E} \cap \{\hbox{$\bar{X}$ visits all sites
of $\lin 0,N\rin$  i.o.}\}.
\end{eqnarray*}
\begin{itemize}
\item[(i)]If  $\mathcal{E}$ (resp. $\mathcal{E}'$) has positive probability, then the VRRW on $\Z$ has positive probability to localize on a
subset of length at most (resp. equal to) $N+1$.
\item[(ii)] Reciprocally, if the VRRW on $\Z$ has  positive probability to localize on a subset of length at most (resp. exactly) $N+1$,
then there exists some initial local time configuration $\mathcal{C}$ such that, for the $\mathcal{C}$-VRRW on
  $\lin 0,N\rin$, the event $\mathcal{E}$ (resp. $\mathcal{E}'$) has positive probability.
\end{itemize}
\end{lem}

\begin{proof}
Define a sequence $(\chi_n)_{n\ge 0}$ of random variables which are,
conditionally on  $\bar{X}$, independent  with law
\begin{equation*}
\pp\{\chi_n=1\ |\ \bar{X}\}=1-\pp\{\chi_n=0\ |\  \bar{X}
\}=\frac{w(0)1_{\{\bar{X}_n=0\}}}{w(0)+w(\bar{Z}_n(1))}+\frac{w(0)1_{\{\bar{X}_n=N\}}}{w(0)+w(\bar{Z}_n(N-1))}.
\end{equation*}
The Borel-Cantelli Lemma applied to the sequence $(\chi_{n})$
yields
\begin{equation}\label{incl}
\mathcal{E}\subset\left\{\sum \chi_n<\infty\right\}.
\end{equation}
Let us also note
that if $\pp\{\sum \chi_n<\infty\}>0$, then necessarily $\pp\{\sum \chi_n=0\}>0$
since we just need to change a finite number of $\chi_n$. Moreover,
it is clear that we can construct a VRRW $X$ on $\Z$ and a reflected random walk $\bar{X}$ on $\lin 0, N\rin$ on the same probability space
in such way that $X$ and $\bar{X}$ coincide on the event $\{\sum \chi_n=0\}$.
Therefore, if $\pp\{\mathcal{E}\}>0$, it follows from \eqref{incl}
that the VRRW $X$ on $\Z$ localizes on $\lin 0,N\rin$ with positive probability.
Moreover, if $\pp\{\mathcal{E}'\}>0$, we find that
$$
\pp\left\{\{\sum \chi_n=0\}\cap\{\hbox{$\bar{X}$ visits all sites of
$\lin 0,N\rin$ i.o.}\}\right\} >0
$$
which implies that, with positive probability, $X$ visits every site of $\lin 0,N\rin$ i.o. without ever exiting the interval.

Reciprocally, if the VRRW on $\Z$ has positive probability to
localize on some interval $\lin x,x+N\rin$, then, clearly, there exists some
initial local time configuration $\mathcal{C}$ on $\Z$ such that the
$\mathcal{C}$-VRRW on $\Z$ has positive probability never to exit
the interval $\lin 0,N\rin$.  On this event, the $\mathcal{C}$-VRRW
on $\Z$ and the restricted $\mathcal{C}$-VRRW on $\lin 0,N\rin$
coincide. We conclude the proof using Lemma \ref{Yfini2} which
implies that, on this event, $Y^+_\infty(0)$ and $Y^{-}_\infty(N)$
are both finite.
\end{proof}

\section{Proof of Theorem \ref{4probapositive}}\label{section4}
We split the proof of the theorem into several propositions. We
start with two elementary observations:

\begin{prop}\label{proptwosite} Let $w$ be a weight sequence. 
\begin{itemize}
\item If $\sum 1/w(k) = \infty$, then we have, a.s., $|R'|\neq
2$.
\item Conversely, if the sum above is finite and the weight sequence $w$ is non-decreasing, then, a.s., $|R'|=
2$.
\end{itemize}
\end{prop}

\begin{proof} Assume that $\sum 1/w(k) = \infty$ and
consider the reflected VRRW
$\bar{X}$ on $\lin 0,1\rin$ starting from some initial configuration
$\mathcal{C}$. We have
$$\bar{Y}_\infty^+(0)=\sum_{k\ge 0} \frac{1_{\{\bar{X}_{k+1}=1\}}}{w(\bar{Z}_k(1))}=\sum_{i=\bar{Z}_0(1)}^\infty \frac{1}{w(i)}=\infty.$$
Thus,  Lemma \ref{lemrest2} implies that $\pp\{|R'| = 2\}=0$.

Reciprocally, if $w$ is non-decreasing and $\sum 1/w(k) < \infty$, then $\sum 1/w(k)^2 <\infty$ and we can invoke Lemma \ref{Yfini} to conclude.
\end{proof}

\begin{prop}
\label{R'3} For any weight sequence $w$, we have, a.s., $|R'|\neq
3$.
\end{prop}
\begin{proof}Consider the reflected VRRW
$\bar{X}$ on $\lin 0,2\rin$ starting from some initial configuration
$\mathcal{C}$. We distinguish two cases:
\begin{itemize}
\item if $\sum 1/w(n)<\infty$, then   Rubin's construction at site 1 implies that either  $0$ or $2$ is visited only finitely many times (notice that we do not require here $w$ to be monotonic).
\item if $\sum 1/w(n)=\infty$, then we have
$$\bar{Y}_\infty^+(0)+\bar{Y}_\infty^-(2)=\sum_{k\ge 0} \frac{1_{\{\bar{X}_{k+1}=1\}}}{w(\bar{Z}_k(1))}=\sum_{i=\bar{Z}_0(1)}^\infty \frac{1}{w(i)}=\infty.$$
\end{itemize}
Thus, in both cases,  Lemma \ref{lemrest2} implies that $\pp\{|R'|= 3\}=0$.

\end{proof}

\begin{rem} \emph{Let us note that the result above does not hold for edge-reinforced random walks: if, for
instance, $w(n)=1$ when $n$ is even and $w(n)=n^2$ when $n$ is odd,
then $R'=\{-1,0,1\}$ a.s., see Sellke \cite{Sel}.}
\end{rem}

In the rest of the paper, given two sequences $(u_n)_{n\ge 1}$ and
$(v_n)_{n\ge 1}$, we shall use Tarr\`es's notation \cite{T1,T2} and
write $u_n\equiv v_n$, when $(u_n-v_n)_{n\ge 1}$ is a converging
sequence.

\begin{lem}\label{lem4_0}
Assume that $w$ is non-decreasing. Then
\begin{eqnarray*}
|R'|=4 \textrm{ with positive probability} \quad \Longrightarrow \quad \sum_{n\ge 0} \frac 1 {w(n)^2} < \infty.
\end{eqnarray*}
\end{lem}
\begin{proof}
Assume that $\sum 1/w(n)^2=\infty$. Let $\bar{X}$ be a VRRW on $\lin
0,3\rin$ starting from some initial configuration $\mathcal{C}$.
Recall that (\ref{eqW}) states that
\begin{eqnarray}
\nonumber W(\bar{Z}_{n}(2))-W(\bar{Z}_{n}(0))&=&\bar{Y}_n^-(3)-\bar{Y}_n^+(-1)+\bar{M}_{n}(1)+ c\\
\label{eqw2}& = & \bar{Y}_n^-(3) +\bar{M}_{n}(1) + c,
\end{eqnarray}
where $c$ is some constant depending on the initial configuration
$\mathcal{C}$. Assume now that $\bar{Y}_\infty^-(3)$ is finite and
let us prove that necessarily $\bar{Y}_\infty^+(0) = \infty$.
Equation \eqref{eqw2} becomes
\begin{equation}\label{eqW3}
W(\bar{Z}_{n}(2))-W(\bar{Z}_{n}(0))\equiv \bar{M}_{n}(1).
\end{equation}
According to Theorem $2.14$ of \cite{HH}, either $\bar{M}_n(1)$
converges or $\limsup \bar{M}_n(1) =-\liminf \bar{M}_n(1) =\infty$.
On one hand, remark that
$$\sum_{n\ge 0}(\bar{M}_{n+1}(1)-\bar{M}_n(1))^2\ge \sum_{n\ge
0}\frac{1_{\{\bar{X}_n=0\}}}{w(\bar{Z}_n(0))^2}.$$ Hence, on the
event $\{\bar{Z}_\infty(0)=\infty\}$, we have $\limsup \bar{M}_n
=-\liminf \bar{M}_n =\infty$. On the other hand, by periodicity,
$\bar{Z}_\infty(2) \vee \bar{Z}_\infty(0) = \infty$. Recalling that
$\lim_{x\rightarrow \infty} W(x)=\infty$, we deduce, using
\eqref{eqW3} that $\{\bar{M}_n(1)\hbox{ converges}\}\subset
\{\bar{Z}_\infty(2)=\infty\}\cap \{\bar{Z}_\infty(0)=\infty\}$.
Therefore, a.s.,
$$
\liminf_n \bar{M}_n(1) = -\infty \quad\hbox{and}\quad \limsup_n
\bar{M}_n(1) = +\infty.
$$
In particular,  there exists a.s. arbitrarily large integers $n$,
such that $W(\bar{Z}_n(0)) \ge W(\bar{Z}_n(2))+1$. Pick such an $n$
and let $m$ be the largest integer smaller than $n$ such that
$W(\bar{Z}_m(0)) \le W(\bar{Z}_m(2))$. For $k\in (m,n]$, we have
$\bar{Z}_k(0)\ge \bar{Z}_k(2)$ hence $$\bar{Z}_k(1)\le
\bar{Z}_k(0)+\bar{Z}_k(2)+\bar{Z}_0(1)\le 3\bar{Z}_k(0),$$ assuming
that $n$ is large enough. Since $w$ is non-decreasing, we get
\begin{eqnarray*}
\sum_{k\in (m,n]} \frac{1_{\{\bar{X}_k=0\}}}{w(\bar{Z}_k(1))} \ge  \sum_{k\in (m,n]} \frac{1_{\{\bar{X}_k=0\}}}{w(3\bar{Z}_k(0))} & = & \sum_{i = \bar{Z}_{m}(0)+1}^{\bar{Z}_{n}(0)}\frac{1}{w(3i)}\\
&\ge & \frac{1}{3} \Big\{W(\bar{Z}_n(0))-W(\bar{Z}_m(0)+1)\Big\}\ge
\frac{1}{4}.
\end{eqnarray*}
As this holds for infinitely many $n$, we deduce that, a.s.
$$\bar{Y}_\infty^{+}(0)=\sum_k \frac{1_{\{\bar{X}_k=0\}}}{w(\bar{Z}_k(1))} = \infty.$$
We conclude by using Lemma \ref{lemrest2}.
\end{proof}

Let us note that Lemma \ref{lem4_0} together with Lemma
\ref{Halphaw} of the appendix imply that, for any non-decreasing
weight sequence $w$, we have
$$
\pp\{|R'| = 4\} >0  \hbox{ or } \alpha_c(w) < \infty \
\Longrightarrow \  \sum_{n}\frac{1}{w(n)^2} < \infty.
$$
Hence, when proving Theorem \ref{4probapositive}, we can assume,
without loss of generality, that $\sum_n 1/w(n)^2 < \infty$. In
particular, the martingales introduced in Section \ref{section3.2}
converge a.s. and in $L^2$.

\begin{prop} \label{prop4pos}
Assume that $w$ is non-decreasing and that \eqref{twosite} holds. We
have
\begin{equation*}
|R'|=4 \textrm{ with positive probability}\ \Longleftrightarrow \ \alpha_c(w)<\infty.\\
\end{equation*}
\end{prop}

\begin{proof} Let us first suppose that $|R'|= 4$ with positive probability.
Thus, according to Lemma \ref{lemrest2}, there exists some initial
local time configuration $\mathcal{C}$ such that, for the
$\mathcal{C}$-VRRW $\bar{X}$ on $\lin 0,3\rin$, the event
$\mathcal{E}:=\{\bar{Y}^+_\infty(0)<\infty\}\cap\{\bar{Y}^{-}_\infty(3)<\infty\}$
has positive probability. Using \eqref{eqW}, we find that
\begin{eqnarray*}
W(\bar{Z}_{n}(2))-W(\bar{Z}_{n}(0))=\bar{Y}_n^-(3)-\bar{Y}_n^+(-1)+\bar{M}_{n}(1)+C\\
W(\bar{Z}_{n}(3))-W(\bar{Z}_{n}(1))=\bar{Y}_n^-(4)-\bar{Y}_n^+(0)+\bar{M}_{n}(2)+C'.
\end{eqnarray*}
As we already noticed, we can assume without loss of generality that
$\sum 1/w(n)^2 < \infty$ so the martingales $\bar{M}_n(1)$ and
$\bar{M}_n(2)$ converge. Hence, there exist finite random variables
$\alpha,\beta$, such that, on the event $\mathcal{E}$,
\begin{eqnarray}\label{eqalphabeta}
W(\bar{Z}_n(1)) - W(\bar{Z}_n(3)) =  \alpha + o(1),\\
\notag W(\bar{Z}_n(2)) - W(\bar{Z}_n(0)) = \beta + o(1).
\end{eqnarray}
For $n$ large enough, this yields
$$\max(\bar{Z}_n(1),\bar{Z}_n(2))\le W^{-1}(W(\max(\bar{Z}_n(0),\bar{Z}_n(3)))+\gamma)$$
with $\gamma:=|\alpha|+|\beta|+1$. Hence, we have
\begin{eqnarray*}
\bar{Y}^+_\infty(0)+\bar{Y}^{-}_\infty(3)&=&\sum_{k\ge 0} \left(\frac {1_{\{\bar{X}_k=0\}}}{w(\bar{Z}_{k}(1))} +  \frac {1_{\{\bar{X}_k=3\}}}{w(\bar{Z}_{k}(2))}\right)\\
&\ge & \sum_{k\ge 0} \frac {1_{\{\bar{X}_k\in \{0,3\}\}}}{w(\max(\bar{Z}_{k}(1),\bar{Z}_{k}(2)))} \\
&\ge & c\sum_{k\ge 0} \frac {1_{\{\bar{X}_k\in \{0,3\}\}}}{w(W^{-1}(W(\max(\bar{Z}_k(0),\bar{Z}_k(3)))+\gamma))}\\
&\ge & c \sum_{k\ge 0} \frac {1}{w(W^{-1}(W(k)+\gamma))}.
\end{eqnarray*}
Therefore, on the event $\mathcal{E}$, we have
$I_{\gamma}(w)<\infty$. This shows that $\alpha_c(w)<\infty$.

\medskip

We now prove the converse implication. Let us assume that
$I_\delta(w)<\infty$ for some $\delta>0$. In particular, $\sum_n
1/w(n)^2<\infty$ (\emph{c.f.} Lemma \ref{Halphaw}). In view of Lemma
\ref{lemrest2}, we will show that, for the reflected VRRW $\bar{X}$
on $\lin 0,3 \rin$, the event
$\{\bar{Y}^+_\infty(0)<\infty\}\cap\{\bar{Y}^{-}_\infty(3)<\infty\}$
has positive probability. This will insure that the VRRW on $\Z$
localizes with positive probability on a subset of size less or
equal to $4$ which will complete the proof of the proposition since
localization on $2$ or $3$ sites is not possible with our
assumptions on $w$.

We use the time-line construction. As explained in the previous
section we can construct $\bar{X}$ from a sequence $(\xi_n^\pm(y),\,
n\ge 0,\, y\in\{1,2\})$ of independent exponential random variables
with mean $1$. Observe that the sequences $(\xi_n^\pm(1)/w(n),n\ge
0)$ define a $w$-urn process via Rubin's construction (choosing $+$
for the red balls). Let $\hat{M}_\infty(1)$ denote the limit of this
urn defined as in Proposition \ref{propurne}. Then, we have
$\hat{M}_\infty(1)\ge \delta+1$ with positive probability. Recall
the definition of $Y^R$ given in Corollary \ref{Yurn} and note that,
on the event $\{\hat{M}_\infty(1)\ge \delta+1\}$, the random
variable $Y^R$ is finite. Besides, using Lemma \ref{couplrest} to
compare $\bar{X}$ with the walk restricted on $\lin 0,2\rin$ (which
correspond to the urn process above), we get
\begin{eqnarray*}
\bar{Y}^+_\infty(0)\le Y^R<\infty \quad \mbox{on the event }
\{\hat{M}_\infty(1)\ge \delta+1\}.
\end{eqnarray*}
By symmetry, considering the limit $\hat{M}_\infty(2)$ of the urn
process $(\xi_n^\pm(2)/w(n),n\ge 0)$, we also find that
\begin{eqnarray*}
\bar{Y}^-_\infty(3) < \infty \quad \mbox{on the event }
\{\hat{M}_\infty(2)\leq  -\delta-1\}.
\end{eqnarray*}
The random variables $\hat{M}_\infty(1)$ and $\hat{M}_\infty(2)$
being independent, we conclude that
$\{\bar{Y}^+_\infty(0)<\infty\}\cap\{\bar{Y}^{-}_\infty(3)<\infty\}$
has positive probability.
\end{proof}

\begin{prop}\label{prop5ex} Assume that $w$ is non-decreasing and that
\eqref{twosite} holds. We have
\begin{equation*}
\alpha_c(w)\in (0,\infty) \ \Longrightarrow \ |R'|=5 \textrm{ with positive probability}.\\
\end{equation*}
\end{prop}

\begin{proof} Assume that $\alpha_c(w)\in (0,\infty)$. In
particular, we have $\sum_{n}1/w(n)^2 < \infty$. Since the walks
associated with a weight $w$ and any non-zero multiple of $w$ have
the same law, we will assume without loss of generality that
$w(0)\ge 1$. Let $\bar{X}$ denote the VRRW reflected on $\lin
0,4\rin$. Let us prove that, with positive probability,
$\bar{Y}^+_\infty(0)$ and $\bar{Y}^-_\infty(4)$ are both finite and
$\bar{X}$ visits all sites of $\lin 0,4 \rin$ infinitely often.

We use again the time-line representation explained in Section
\ref{section3} except that we will change the construction slightly
for the transition at site $2$. Recall that according to the
original construction, when the process jumps for the $k$-th time
from $1$ (resp. $3$) to $2$ and has made $m$ visits to $1$ (resp.
$3$) before time $t$, then we attached to the oriented edge $(2,1)$
(resp. $(2,3)$) a clock which rings after time $\xi_k^-(2)/w(m)$
(resp. $\xi_k^+(2)/w(m)$). In our new construction, we choose to
attach instead a clock which rings after time $\xi_m^-(2)/w(m)$
(resp. $\xi_m^+(2)/w(m)$). The random variables $(\xi_k^\pm(2),\,
k\geq 0)$ being i.i.d, this modification does not change the law of
$\bar{X}$ (some random variables $\xi_k^\pm(2)$ are simply never
used).

Fix some $0<\varepsilon < 1$ and  consider the two $w$-urn processes
$u_1:=(\xi_n^\pm(1)/w(n),n\ge 0)$, and $u_3:=(\xi_n^\pm(3)/w(n),n\ge
0)$. Since $\bar{Y}^+_\infty(0)$ is stochastically smaller than it
would be for the process reflected in $\lin 0,2\rin$, using  similar
arguments as in the proof of Proposition \ref{prop4pos}, we see that
there exists a set $E_1\subset (\R_+^2)^\N$, such that the event
$\mathcal{E}_1:=\{u_1\in E_1\}$ has positive probability and on
which $\bar{Y}^+_\infty(0)\le \varepsilon^3$. By symmetry, there
exists a set $E_2$, such that $\mathcal{E}_2:=\{u_3\in E_2\}$ has
positive probability and on which $\bar{Y}^-_\infty(4)\le
\varepsilon^3$. By independence of the urns $u_1$ and $u_3$, the
event $\mathcal{E}_1\cap \mathcal{E}_2$ also has positive
probability. In view of Lemma \ref{lemrest2}, it remains to prove
that, on this event, $\bar{X}$ visits all the sites of $\lin 0,4
\rin$ infinitely often with positive probability.

We now consider the urn process $u_2:=(\xi_n^\pm(2)/w(n),n\ge 0)$.
Recall that, according to \eqref{limitUrn}, we may express the limit
$\hat{M}_\infty(2)$ of this urn in the form:
\begin{equation}\label{eqs1}
\hat{M}_\infty(2)= \sum_{n \ge 0} \left( \frac {1-\xi_n^+(2)}{w(n)}
\right) - \sum_{n\ge 0} \left(\frac{1-\xi_n^-(2)}{w(n)}\right).
\end{equation}
Similarly, it is not difficult to check that we can
also express the limit of the martingale $\bar{M}_\infty(2) :=
\lim_{n\to \infty} (\bar{Y}^+_n(2)-\bar{Y}^-_n(2))$ in the form
\begin{equation}\label{eqs2}
\bar{M}_\infty(2)= \sum_{n \ge 0} \left( \frac
{1-\xi_{c_n}^+(2)}{w(c_n)} \right) - \sum_{n\ge 0}
\left(\frac{1-\xi_{d_n}^-(2)}{w(d_n)}\right).
\end{equation}
where $(c_n,n\ge 0)$ and $(d_n,n\ge 0)$ are the increasing (random)
sequences such that $\bar{Y}^+_n(2)= \sum_{c_k \leq n } 1/w(c_k)$
and $\bar{Y}^-_{n}(2)=\sum_{d_k\leq n} 1/w(d_k)$. The idea now is to
compare $\bar{M}_\infty(2)$ and $\hat{M}_\infty(2)$ and prove that,
on the event $\mathcal{E}_1\cap \mathcal{E}_2$ their values are
close. Then we will use the fact that $\hat{M}_\infty(2)$ has a
density to deduce that $\bar{M}_\infty(2)$ can be smaller than
$\alpha_c(w)$.

Subtracting \eqref{eqs2} from \eqref{eqs1}, we find that
$$
\hat{M}_\infty(2) - \bar{M}_\infty(2) = \sum_{n\ge 0} \left(\frac
{1-\xi_{i_n}^+(2)}{w(i_n)}\right) + \sum_{n\ge 0}
\left(\frac{1-\xi_{j_n}^-(2)}{w(j_n)}\right),
$$
where $(i_n,n\ge 0)$ and $(j_n,n\ge 0)$ are the complementary
sequences of $(c_n,n\ge 0)$ and $(d_n,n\ge 0)$. Moreover, using
relation \eqref{Ynpm}, we have
$$
\bar{Y}^+_\infty(0) = \sum_n \frac{1}{w(j_n)} \quad\hbox{ and }\quad
\bar{Y}^-_\infty(4) = \sum_n \frac{1}{w(i_n)}.
$$
But, using similar arguments as in the proof of \eqref{Tx}, we obtain
\begin{multline*}
\E\left[\left(\sum_{n\ge 0} \frac
{1-\xi_{i_n}^+(2)}{w(i_n)}\right)^2\ \Big|\ \mathcal{E}_1\cap
\mathcal{E}_2\right] =  \E\left[ \sum_{n\ge 0} \frac 1 {w(i_n)^2}\
\Big|\ \mathcal{E}_1\cap \mathcal{E}_2\right] \\
\le \frac{1}{w(0)}\E\left[\sum_{n\ge 0} \frac 1 {w(i_n)}\ \Big|\
\mathcal{E}_1\cap \mathcal{E}_2\right]  \le \E[\bar{Y}_\infty^-(4)\mid
\mathcal{E}_1\cap \mathcal{E}_2] \le \varepsilon^3.
\end{multline*}
Using Tchebychev's inequality, we deduce
$$\pp\left\{|\hat{M}_\infty(2)-\bar{M}_\infty(2)|\ge 2\varepsilon\mid \mathcal{E}_1\cap \mathcal{E}_2 \right\} \le \varepsilon. $$
Recalling that $\hat{M}_\infty(2)$ has a density with support on the
whole of $\R$ (\emph{c.f.} Proposition \ref{propurne}), we can pick
$\eta>0$ such that $\pp\{|\hat{M}_\infty(2)|\le \eta\} =
2\varepsilon$. This yields
$$\pp\Big\{|\bar{M}_\infty(2)|\ge \eta + 2\varepsilon \mid \mathcal{E}_1\cap \mathcal{E}_2\Big\} \le 1-2\varepsilon + \varepsilon \le 1-\varepsilon,$$
so the set $\mathcal{E}_3:=\mathcal{E}_1\cap \mathcal{E}_2\cap
\{|\bar{M}_\infty(2))| \le \eta + 2\epsilon\}$ has positive probability.
Moreover, we have,
$$W(\bar{Z}_n(1))-W(\bar{Z}_n(3))= \bar{Y}^+_n(0) - \bar{Y}^-_n(4) -\bar{M}_n(2),$$
and therefore, for all $n$ large enough, on $\mathcal{E}_3$,
\begin{equation*}
\bar{Z}_n(2) \le  \bar{Z}_n(1) + \bar{Z}_n(3)  \le  \bar{Z}_n(1) +
W^{-1}(W(\bar{Z}_n(1))+\eta +4\varepsilon) .
\end{equation*}
Choosing $\varepsilon$ small enough such that $\delta:=\eta
+4\varepsilon <\alpha_c(w)$,  we obtain, for $N$ large
\begin{eqnarray*}
\bar{Y}^+_\infty(1)&\ge& \sum_{n\ge N}
\frac{1_{\{\bar{X}_n=1,\bar{X}_{n+1}=2\}}}{w(\bar{Z}_n(1) +
W^{-1}(W(\bar{Z}_n(1))+\delta))}\\
& \ge &  \sum_{n\ge N} \frac{1}{w(n +
W^{-1}(W(n)+\delta))}-\sum_{n\ge N}
\frac{1_{\{\bar{X}_n=1,\bar{X}_{n+1}=0\}}}{w(\bar{Z}_n(1))}\\
& \ge &  \sum_{n\ge N} \frac{1}{w(n +
W^{-1}(W(n)+\delta))}-\bar{Y}^+_\infty(0).\\
\end{eqnarray*}
It follows from Lemma \ref{teclem} of the Appendix that
$\bar{Y}^+_\infty(1)$ is infinite on $\mathcal{E}_3$. By symmetry,
we also have $\bar{Y}^-_\infty(3)=\infty$ on $\mathcal{E}_3$.
Therefore, according to Lemma \ref{Yfini}, on $\mathcal{E}_3$, the
VRRW on $\lin 0,4\rin$ visits every site infinitely often. This
concludes the proof of the proposition.
\end{proof}

\begin{prop}\label{Prop45as}
Assume that $w$ is non-decreasing and that \eqref{twosite} holds.
Then
\begin{eqnarray*}
\alpha_c(w)<\infty \quad \Longrightarrow \quad |R'|\in \{4,5\} \textrm{ almost surely.}
\end{eqnarray*}
\end{prop}
\begin{proof}
We first argue that, if $\alpha_c(w)<\infty$, then $R'$ is a.s.
finite and non-empty. Indeed, recalling Lemma \ref{couplrest}, each
time $X$ visits a new site, say $x>0$, as long as it does not visit
$x+2$, the restriction of $X$ to the set $\{x-1,x,x+1\}$ can be
coupled with a $w$-urn process in such a way that it always makes
less jumps to $x+1$ than the urn process. Then, Corollary \ref{Yurn}
and Lemma \ref{Yfini} insure that $X$ never visits $x+2$ with a
positive probability uniformly bounded from below by a constant
depending only on this urn process (and therefore which does not
depend on the past trajectory of $X$ before its first visit to $x$).
It follows that, a.s., $\limsup X_n < \infty$ and by symmetry
$\liminf X_n > -\infty$. Hence the walk localizes on a finite set
almost surely.

\medskip

Let us now assume, by contradiction, that $|R'|=N+1\ge 6$ with
positive probability. Thus, according to Lemma \ref{lemrest2}, there
exists some initial local time configuration $\mathcal{C}$ such
that, for the $\mathcal{C}$-VRRW  $\bar{X}$ on $\lin 0,N\rin$, the
event
$$\mathcal{E}:=\{\bar{Y}_\infty^+(0)+\bar{Y}_\infty^-(N)<\infty\}\cap\{\bar{X}
\mbox{ visits 0 and $N$ i.o.} \}$$ has positive probability.
Moreover, using \eqref{eqW}, we have
\begin{eqnarray*}
W(\bar{Z}_{n}(2))-W(\bar{Z}_{n}(0))&\equiv &\bar{Y}_n^-(3)\\
W(\bar{Z}_{n}(3))-W(\bar{Z}_{n}(1))&\equiv &\bar{Y}_n^-(4)-\bar{Y}_n^+(0).
\end{eqnarray*}
On the event $\mathcal{E}$, the quantity $\bar{Y}_\infty^+(0)$ is
finite whereas $\bar{Y}_\infty^-(3)$ and $\bar{Y}_\infty^-(4)$ are
infinite according to Lemma \ref{Yfini} since $N\geq 5$. Thus, for
all $A>0$, the stopping time
\begin{eqnarray*}
T_A:=\inf \left\{n \ge 0 \ : \ \begin{array}{l}  \bar{X}_n=2\\
                            W( \bar{Z}_n(0)) \le W( \bar{Z}_n(2)) - A  \\
                               W( \bar{Z}_n(1)) \le W( \bar{Z}_n(3)) - A
                              \end{array}
\right\}
\end{eqnarray*}
is finite on $\mathcal{E}$. We claim that
\begin{eqnarray}\label{Yinftyfini}
\pp\{ \bar{Y}_\infty^+(1)=\infty\mid T_A<\infty \} \to 0 \quad \textrm{as }A\to \infty.
\end{eqnarray}
For the time being, assume that \eqref{Yinftyfini} holds. As before,
Lemma \ref{Yfini} states that, on $\mathcal{E}$, we have
$\bar{Y}_\infty^+(1) = \infty $. Thus, for all $A>0$, we get
\begin{equation*}
\pp\{\mathcal{E}\}\le \pp\left\{\{\bar{Y}_\infty^+(1)=\infty\} \cap
\{T_A<\infty\}\right\}\le  \pp\{ \bar{Y}_\infty^+(1)=\infty\mid
T_A<\infty \},
\end{equation*}
which yields $\pp\{\mathcal{E}\}=0$ and contradicts the initial
assumption that the walk localizes with positive probability on
 more than $5$ sites.

\medskip

It remains to prove \eqref{Yinftyfini}. For $A>0$, consider a
process $\bar{X}^A$ which is, up to time $T_A$, equal to the VRRW
$\bar{X}$ on $\lin 0,N\rin$ and which, after time $T_A$, has the
transition of the VRRW restricted on $\lin 0,2\rin$. In view of
Lemma \ref{couplrest}, we can construct $\bar{X}^A$ together with
$\bar{X}$  in such way that, with obvious notation,
$$\bar{Y}_\infty^+(0)\le \bar{Y}^{A,+}_\infty(0) \quad \mbox{ and }\quad \bar{Y}_\infty^+(1)\le \bar{Y}^{A,+}_\infty(1).$$
After time $T_A$, the process $\bar{X}^A$ is simply an urn process.
Hence, using the same arguments as in the proof of Proposition
\ref{propurne}, for $n\geq T_A$, the process
$$\hat{M}^A_n := W( \bar{Z}^A_n(0))- W( \bar{Z}^A_n(2))$$
is a martingale with quadratic variation bounded by $2\sum_n
1/w(n)^2$. Noticing that, by definition of $T_A$, we have
$\hat{M}^A_{T_A}\le -A$, the maximal inequality for martingales
shows that, for any $\varepsilon
>0$, there exists a constant $C>0$, such that for all $A>0$,
\begin{eqnarray}
\label{TA} \pp\left\{\sup_{n\ge T_A}\, W(\bar{Z}^A_n(0))- W(
\bar{Z}^A_n(2))\ge -A+ C \;|\; \mathcal{F}_{T_A}\right\} \le
\varepsilon.
\end{eqnarray}
Moreover, for every odd integer $n \ge T_A$, we have $\bar{X}^A_n =
1$ from which we deduce that, for all $n \ge T_A$,
$$\bar{Z}^A_n(1) \;\ge\; \bar{Z}^A_n(2)+\bar{Z}^A_n(0)+\bar{Z}^A_{T_A}(1)-\bar{Z}^A_{T_A}(0)-\bar{Z}^A_{T_A}(2)\;\ge\; \bar{Z}^A_n(2)-\bar{Z}^A_{T_A}(2).$$
On the event $\{\sup_{n\ge T_A}\, W(\bar{Z}^A_n(0))- W(
\bar{Z}^A_n(2))< -A+ C\}$, we get, for $n\ge T_A$,
$$\bar{Z}^A_n(1)\ge W^{-1}(W(\bar{Z}^A_n(0)+A-C) -\bar{Z}^A_{T_A}(2).$$
This yields
\begin{eqnarray*} \bar{Y}^{A,+}_\infty(0) &=&\bar{Y}^{A,+}_{T_A}(0) + \sum_{n > T_A}
\frac{1_{\{\bar{X}_n^A = 0\}}}{w(\bar{Z}_n^A(1))}\\
&\leq&\bar{Y}^{A,+}_{T_A}(0) +\sum_{n\ge 0}
\frac{1}{w(W^{-1}(W(n)+A-C)-\bar{Z}^A_{T_A}(2))}
\end{eqnarray*}
with the convention $w(x)=w(0)$ for $x\le 0$. Thus, according to
Lemma \ref{lemaa} of the appendix, for $A >\alpha_c(w)+ C$, we have
$\bar{Y}^{A,+}_\infty(0)<\infty$ on the event
$\{T_A<\infty\}\cap\{\sup_{n\ge T_A}\, W(\bar{Z}^A_n(0))- W(
Z^A_n(2))< -A+ C\}$. Using \eqref{TA}, we obtain
\begin{eqnarray}\label{eka}
\pp\left\{\bar{Y}^{A,+}_\infty(0)=\infty\,|\, T_A<\infty\right\}\le
\varepsilon.
\end{eqnarray}
We can now choose $A_0 > \alpha_c(w) + C$ and $K>0$ such that
$$\pp\{\bar{Y}^{A_0,+}_\infty(0)\ge K\, | \, T_{A_0}<\infty \} \le 2\varepsilon.$$
Notice that for $A>A_0$, the random variable
$\bar{Y}^{A,+}_\infty(0)$ is stochastically dominated by
$\bar{Y}^{A_0,+}_\infty(0)$ (this is again a consequence of Lemma
\ref{couplrest} using the same time-line construction for
$\bar{X}^A$ and $\bar{X}^{A_0}$). Moreover, by hypothesis,
$$\pp\{T_A<\infty \}\ge \pp\{\mathcal{E}\}:= c >0,$$
hence
\begin{equation}
\label{Y((A))+}
\forall A>A_0, \quad \pp\{\bar{Y}^{A,+}_\infty(0)\ge K\, |\,  T_{A}<\infty \} \le 2\varepsilon/c.
\end{equation}

Finally, we consider a third process $\tilde{X}^A$ which coincides
up to time $T_A$ with $\bar{X}$ and $\bar{X}^A$, and which, after
time $T_A$, has the transition of the VRRW restricted on $\lin
0,3\rin$. Again, we can construct these processes in such way that
 $\bar{Y}_\infty^{A,+}(0)$ stochastically dominates $\tilde{Y}_\infty^{A,+}(0)$.
This domination implies, that \eqref{Y((A))+} also holds with
$\tilde{Y}^{A,+}_\infty(0)$ in place of $\bar{Y}^{A,+}_\infty(0)$.
Moreover,
$$\tilde{M}^A_n(2):=W(\tilde{Z}^A_n(3))-W(\tilde{Z}^A_n(1))+\tilde{Y}^{A,+}_n(0) \qquad n\ge T_A,$$
is a martingale with bounded quadratic variation. As before, we
deduce from the maximal inequality for martingales, that, for some
constant $C'>0$ depending only on $\varepsilon$ and the weight
function $w$,
$$\pp\left\{\inf_{n\ge T_A}\, \tilde{M}^A_n(2)-\tilde{M}^A_{T_A}(2)\le -C'\ \Big| \ \kF_{T_A}\right\} \le \varepsilon.$$
Using the facts that $\tilde{M}^A_{T_A}(2)\ge A$ and that
$\tilde{Y}^{A,+}_n(0)\ge K$ with probability smaller than
$2\varepsilon/c$ on the event $\{T_{A}<\infty \}$, we obtain
$$\pp\left\{\inf_{n\ge T_A}\, W(\tilde{Z}^A_n(3))-W(\tilde{Z}^A_n(1))\le -C'-K+A\ \Big| \ T_A<\infty\right\} \le \varepsilon',$$
with $\varepsilon'=\varepsilon(1+2/c)$. We now fix $A$ large enough
such that $A-C'-K>\alpha_c(w)$. Using the trivial relation
$\tilde{Z}_n^{A}(2)\ge \tilde{Z}_n^{A}(3)-\tilde{Z}_{T_A}^{A}(3)$
for $n\ge T_A$, we deduce, in the same way as for the proof of
\eqref{eka}, that
$$\pp\left\{ \tilde{Y}^{A,+}_\infty(1)=\infty \ \Big|\ T_A<\infty \right\} \le \varepsilon'.$$
We conclude the proof of \eqref{Yinftyfini} by noticing that
$\tilde{Y}^{A,+}_\infty(1)$  stochastically dominates
$\bar{Y}_\infty^{+}(1)$.
\end{proof}

\begin{lem}\label{sansatome} Assume that $w$ is non-decreasing and
that $\sum_{n} 1/w(n)^2 < \infty$. Fix $\infty\le a<0<b\le \infty$
and let $\bar{X}$ be a $\mathcal{C}$-VRRW on $\lin a,b\rin$ for some
initial local time configuration $\mathcal{C}$. Set
$$\delta:=\lim_{n\rightarrow\infty} W(\bar{Z}_n(1))-W(\bar{Z}_n(-1)) \quad
\mbox{ when the limit exists.}$$ Then, for any $\delta_0\in \R$, we
have
$$\pp\big\{\{\delta \hbox{ exists and equals }\delta_0\} \cap \{\bar{Z}_\infty(0)=\infty\}\big\}=0.$$
\end{lem}

\begin{proof}  Since the weights $w$ and $\lambda w$ (for $\lambda>0$) define the same VRRW, we assume, without loss of generality
that $w(Z_0(1))=1$. Recalling the time-line construction described in Section \ref{section3}, we create the $\mathcal{C}$-VRRW $\bar{X}$ on $\lin a,b\rin$
from a collection  $((\xi_n^\pm(y),n\ge 0),y\in \ilin a,b\irin)$. Set
$$\kH:=\left((\xi_n^\pm(y),n\ge 0, y\neq 0),(\xi_n^-(0),n\ge 0),(\xi_n^+(0),n\ge 1)\right)\in \R_+^\N$$
and let $\mu$ denote the product measure on $\R_+^\N$ under which
$\kH$ is a collection of i.i.d. exponential random variables with
mean $1$. Then, given $\kH$ and some other variable $\xi_0^+(0)$, the pair $(\kH,\xi_0^+(0))$ defines a process $\bar{X} = \bar{X}(\kH,\xi_0^+(0))$ via the time-line construction which is a VRRW under the product probability $\pp :=\mu\times\hbox{Exp}(1)$.
For $u>0$, define
$$ \mathcal{B}_{u}=\{\kH\in \R_+^\N,\; \delta(\kH,u) \hbox{ exists and equals }\delta_0 \hbox{ and } \bar{Z}_\infty(0)=\infty\}$$
and
$$
\mathcal{B} =\{(\kH,u) \in \R_+^\N\times\R_+,\; \kH \in \mathcal{B}_u\}.
$$
We will prove that, for almost every $u>0$ and $h>0$,
\begin{equation}\label{intersecB}
\mu\{\mathcal{B}_u\cap \mathcal{B}_{u+h}\}=0.
\end{equation}
This equation implies that $\{u\in \R^+,\mu\{\mathcal{B}_u\}>0\}$
has zero Lebesgue measure. Hence
$$\pp\{\{\delta=\delta_0\}\cap\{\bar{Z}_\infty(0)=\infty\}\}= \pp\{\mathcal{B}\} =\int_0^\infty e^{-u}\mu\{\mathcal{B}_u\}du=0.$$
It remains to prove \eqref{intersecB}.
Given  $(\kH,u)$ such that $\bar{Z}_\infty(1)= Z_\infty(-1)=\infty$, we can define, as in the proof of Lemma \ref{Yfini}, the increasing sequences
$(i_k^\pm,k\ge 0)$ such that, for all $n\ge 0$,
$$Y_n^\pm(0)= \sum_{i_k^\pm\le n} \frac{1}{w(i_k^\pm)}.$$
 For $t>0$, define also
\begin{equation}\label{defzip}
z^\pm_t:=\inf\left\{n\ :\ \sum_{i_k^\pm \le n} \frac {\xi_k^\pm(0)}{w(i_k^\pm)}>t\right\}.
\end{equation}
Thus, $z^\pm_t$ represents the local time at site $\pm 1$ when the
clock process  attached to  site $0$ has consumed a time $t$. Hence,
another way to define $z^\pm_t$ is to consider the continuous-time
process $(\tilde{X}(s),s\ge 0)$ associated with $(\kH,u)$ via the
time-line construction (recall that $\bar{X}$ is deduced from
$\tilde{X}$ by a change of time). Defining
$$\tau_t:=\inf\Big\{s>0, \int_{0}^s1_{\{\tilde{X}_s=0\}}ds > t\Big\},$$
we get that $z^\pm_t=\tilde{Z}_{\tau_t}(\pm 1)$.

Let us  notice that, on $\mathcal{B}$, $\pp$-a.s., we have $Z_\infty(-1)=Z_\infty(1)=\infty$ so the sequences $(i_k^{\pm})$ are well defined for all $k\ge 0$. Moreover, on the event $\{Z_\infty(-1)=Z_\infty(1)=\infty\}$, the total time consumed by the clock process at site $0$ is infinite
 $\pp$-a.s. (see the proof of Lemma \ref{Yfini}). Hence, on $\mathcal{B}$, the random variables $z^\pm_t$ are finite for all $t>0$, $\pp$-a.s. Define
$$\delta_{t} := W(z^+_t)-W(z^-_t)=W(\tilde{Z}_{\tau_t}(1))-W(\tilde{Z}_{\tau_t}(-1)).$$
By definition of $\delta$, we get
\begin{equation*}
\lim_{t \rightarrow \infty} \delta_t=\delta.\qquad \pp\mbox{-a.s. on  $\mathcal{B}$}.
\end{equation*}
Thus, for almost any $u>0$ (with respect to the Lebesgue measure),
we have
\begin{equation}\label{dtmuas}
\lim_{t \rightarrow \infty} \delta_t(\kH,u)=\delta(\kH,u)\qquad \mbox{for $\mu$-a.e. $\kH\in\mathcal{B}_u$}.
\end{equation}
Now let $u, h > 0$ be fixed and such that \eqref{dtmuas} holds for
$u$ and $u+h$. Pick $\kH \in \mathcal{B}_u\cap\mathcal{B}_{u+h}$.
Lemma \ref{tarres1} implies that, for all $k\ge 0$, we have
$$i_k^+(\kH,u+h)\le i_k^+(\kH,u)  \quad\hbox{ and }\quad i_k^-(\kH,u+h)\ge i_k^-(\kH,u).$$
Recalling that $w(\bar{Z}_0(1))=w(i_0^+)=1$, we deduce from \eqref{defzip} that, for $t>0$,
$$
z^+_t(\kH,u+h)\le z^+_{t-h}(\kH,u) \quad\hbox{ and }\quad z^-_t(\kH,u+h)\ge z^-_t(\kH,u).
$$
This yields
\begin{eqnarray*}
\delta_t(\kH,u)-\delta_t(\kH,u+h)&\ge& W(z^+_{t}(\kH,u))- W(z^+_{t-h}(\kH,u))\\
& \ge & \sum_{z^+_{t-h}\le k< z^+_{t}}\  \frac{1}{w(k)} \\
& \ge & \sum_{z^+_{t-h}\le i_k^+< z^+_{t}}\  \frac{1}{w(i_k^+)} :=
\Delta_{u,h}^+(t),
\end{eqnarray*}
where $z^+_{t-h}, z^+_{t}$ and $i_k^+$ stand for $z^+_{t-h}(\kH,u), z^+_{t}(\kH,u)$ and $i_k^+(\kH,u)$.
In view of \eqref{dtmuas}, we deduce that, for almost every $u,h>0$, we have
$$\mu\{\mathcal{B}_u\cap \mathcal{B}_{u+h}\}
\le \mu\big\{\mathcal{B}_u\cap\{\limsup_{t\to\infty}\Delta_{u,h}^+(t)=0 \}\big\}.$$
It remains to prove that the r.h.s. in the previous inequality is equal to zero. For $\kH\in\mathcal{B}_u$, the quantity
$$h_t^*(\kH,u)  := \sum_{z^+_{t-h}\le i_k^+< z^+_t}\  \frac{\xi^+_k(0)}{w(i_k^+)}$$
is well defined. Moreover, it is clear that,
\begin{equation}\label{zr1}
\lim_{t \rightarrow \infty} h^*_t=h\quad\hbox{  $\mu$-a.s. on
$\mathcal{B}_u$.}
\end{equation}
On the other hand, we have,
$$|h^*_t-\Delta_{u,h}^+(t)|^21_{\mathcal{B}_u}\le 2 \left(\sum_{i_k^+\ge z^+_{t-h}}  \!\!\!\!\frac{1-\xi^+_k(0)}{w(i_k^+)}1_{\{\theta_k<\infty\}}\right)^{\!\!\! 2}+  2\left(\sum_{i_k^+\ge z^+_{t}}\!\!\!\!  \frac{1-\xi^+_k(0)}{w(i_k^+)}1_{\{\theta_k<\infty\}}\right)^{\!\!\! 2},$$
where $\theta_k$ denotes the time of the $k$-th jump of $\bar{X}$
from $1$ to $0$. Let $(\tilde{\mathcal{F}}_t,t>0)$ denote the
natural filtration of the continuous time process
$\tilde{X}_t(\cdot,u)$. Using the same argument as in the proof of
\eqref{Tx}, we find that
\begin{eqnarray*}
\E_\mu\left[(h^*_t-\Delta_{u,h}^+(t))^21_{\mathcal{B}_u}\ \Big| \  \tilde{\mathcal{F}}_{\tau_{t-h}}\right]
&\le &
4\E_\mu\left[\sum_{i_k^+\ge z^+_{t-h}}\  \frac{1_{\{\theta_k<\infty\}}}{w(i_k^+)^2}\ \Big| \ \tilde{\mathcal{F}}_{\tau_{t-h}} \right] \\
&\le & \sum_{k\ge z^+_{t-h}} \frac 4 {w(k)^2}.
\end{eqnarray*}
This yields
$$\mu\left\{|h^*_t-\Delta_{u,h}^+(t)|1_{\mathcal{B}_u}\ge h/2\ \Big| \  \ \tilde{\mathcal{F}}_{\tau_{t-h}}\right\} \le \frac{16}{h^2}\,
 \sum_{k\ge z^+_{t-h}} \frac 1 {w(k)^2}.$$
Hence, by monotone convergence,
\begin{equation}\label{zr2}
\lim_{t\to\infty}\mu\left\{|h^*_t-\Delta_{u,h}^+(t)|1_{\mathcal{B}_u}\ge h/2\right\} = 0
\end{equation}
Combining \eqref{zr1} and \eqref{zr2}, we conclude that
\begin{multline*}
\lim_{t\to\infty}\mu\left\{ \mathcal{B}_u\cap\{\Delta_{u,h}^+(t)
\leq h/4\} \right\} \\ \leq \lim_{t\to\infty}\mu\left\{
\mathcal{B}_u\cap\{|h^*_t-\Delta_{u,h}^+(t)|\geq h/2\} \right\}+
\lim_{t\to\infty}\mu\Big\{ \mathcal{B}_u\cap\{h^*_t\leq 3h/4\}
\Big\} = 0,\end{multline*} which implies
$\mu\big\{\mathcal{B}_u\cap\{\limsup_{t\to\infty}\Delta_{u,h}^+(t)=0
\}\big\} = 0$.
\end{proof}

We can now prove the last part of Theorem \ref{4probapositive}.

\begin{prop}\label{propnot5}
Assume that $w$ is non-decreasing and that \eqref{twosite} holds.
Then
\begin{eqnarray*}
\alpha_c(w)=0 \quad \Longrightarrow \quad |R'|\neq 5 \textrm{ almost surely.}
\end{eqnarray*}
\end{prop}
\begin{proof}
Assume by contradiction that $\alpha_c(w)=0$ and that $|R'|=5$ holds
with positive probability. Thus there exists an initial local time
configuration $\mathcal{C}$ such that, for the $\mathcal{C}$-VRRW
$\bar{X}$ on $\lin 0,4\rin $, the event
$$\mathcal{E}:=\{\bar{Y}_\infty^+(0)+\bar{Y}_\infty^-(4)<\infty\}\cap\{\bar{X}
\mbox{ visits 0 and $4$ i.o.} \}$$ has positive probability.
Moreover, Equation \eqref{eqW}  yields, for $n\ge 0$,
$$W(\bar{Z}_n(1))- W(\bar{Z}_n(3))= \bar{Y}_n^+(0) -\bar{M}_n(2) -\bar{Y}_n^-(4)+c,$$
for some constant $c$ depending on the initial configuration. On the
event $\mathcal{E}$, each term on the r.h.s. of this equation
converges to a limit, thus
 $$\lim_{n\to\infty} W(\bar{Z}_n(1))- W(\bar{Z}_n(3))=\bar{Y}_\infty^+(0) -\bar{M}_\infty(2) -\bar{Y}_\infty^-(4)+c =: \delta$$
 exists and is finite. Moreover, Lemma \ref{sansatome} implies that $\pp\{\{\delta=0\} \cap \mathcal{E}\}=0$.
 Let us now prove that the event $ \mathcal{E}\cap
\{\delta>0\}$ has probability 0  (the same result holds for
$\delta<0$ by symmetry). On this event, for $n$ large enough, we get
\begin{eqnarray}
\label{x1x3} W(\bar{Z}_n(1))\ge W(\bar{Z}_n(3))+\delta',
\end{eqnarray}
with $\delta'=\delta/2$. Besides, we have
\begin{eqnarray*}
 W(\bar{Z}_n(2))-W(\bar{Z}_n(0))\equiv \bar{Y}_n^-(3)\\
 W(\bar{Z}_n(2))-W(\bar{Z}_n(4))\equiv \bar{Y}_n^+(1).
\end{eqnarray*}
Since, on $\mathcal{E}$, the quantities $\bar{Y}_\infty^-(3)$ and
$\bar{Y}_\infty^+(1)$ are infinite (\emph{c.f.} Lemma \ref{Yfini}),
we deduce that $\bar{Z}_n(2)$ is larger than $\bar{Z}_n(0)$ and
$\bar{Z}_n(4)$ for $n$ large enough. Moreover, by periodicity, the
sum of these three quantities is, up to a constant, equal to $n/2$.
Thus, for $n$ large enough, we obtain
$$3\bar{Z}_n(2)\ge \frac{n}{2} \ge \bar{Z}_n(1).$$
Using \eqref{x1x3}, we get, for $n$ large enough, on the event $
\mathcal{E}\cap \{\delta>0\},$
\begin{eqnarray*}
\frac{1_{\{\bar{X}_n=3,\,\bar{X}_{n+1}=2 \}}}{w(\bar{Z}_n(2))} &\le
& \frac{1_{\{\bar{X}_n=3\}}}{w(\bar{Z}_n(1)/3)} \le
\frac{1_{\{\bar{X}_n=3\}}}{w(W^{-1}(W(\bar{Z}_n(3))+\delta' )/3)}.
\end{eqnarray*}
Since $\alpha_c(w)=0$, it follows from Lemma \ref{lemaa} of the
Appendix  that $\bar{Y}_\infty^-(3)$ is finite which contradicts the
fact that $\bar{X}$ visits all the sites of $\lin 0,4 \rin$
infinitely often.
\end{proof}

\medskip
Theorem \ref{4probapositive} is now a consequence of Propositions
\ref{prop4pos}, \ref{prop5ex}, \ref{Prop45as} and \ref{propnot5}.
Let us conclude this section by remarking that we can also describe the
shape of the asymptotic local time configuration when $\alpha_c(w)
<\infty$. Indeed, collecting the results obtained during the proof of
the theorem, it is not difficult to check (the details being left
out for the reader) that the asymptotic local time profile of the
walk on $R'$ at time $n$ takes the form:
\bigskip

\begin{center}
\includegraphics[width=12.6cm,angle=0]{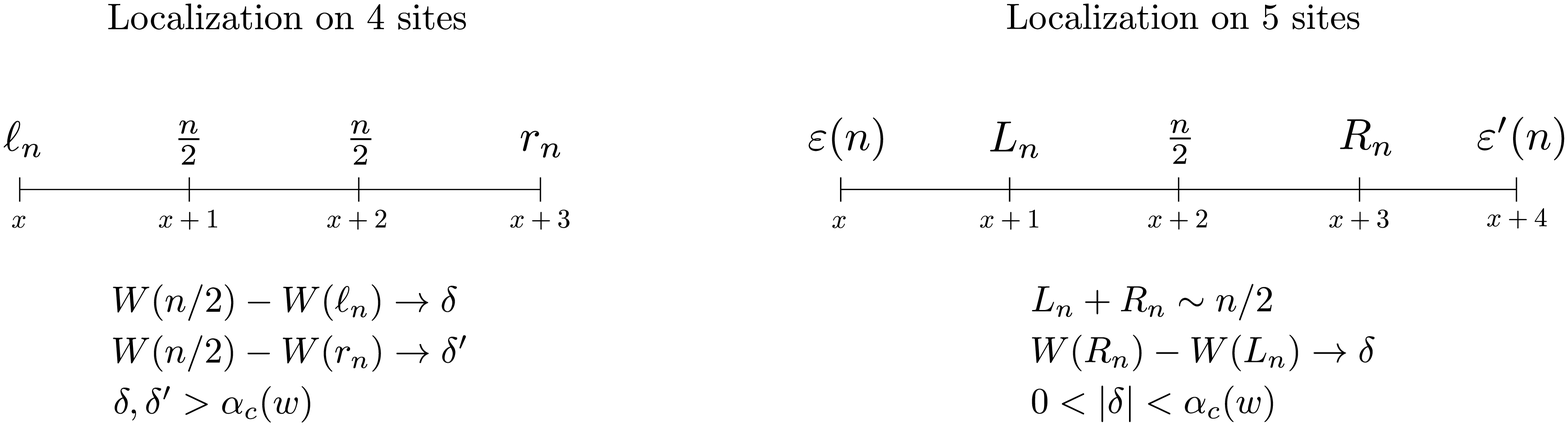}
\end{center}

\bigskip
\noindent In particular, when the walk localizes on $4$ sites, only
the two central sites are visited  a non-negligible proportion of
time (this follows from \eqref{andthelast} of Lemma \ref{lemaa} of
the appendix). When the walk localizes on $5$ sites, a more unusual
behaviour may happen. If the weight function is regularly varying
(for example $w(n)\sim n\log \log n$), then, again, the walk spends
asymptotically all its time on two consecutive vertices:
$$\lim_{n\rightarrow\infty } \frac{L_n\wedge R_n}{n}=0 \qquad \hbox{ and }\qquad  \lim_{n\rightarrow\infty } \frac{L_n\vee R_n}{n}=\frac{1}{2}.$$
However, this result is not true for general weight functions. In
fact, the ratio $Z_n(y)/n$ of time spent at site $y$ may not
converge. For instance, considering the  weight sequence $w_0$ of
Remark \ref{remn!} of the appendix, we find that when the walk
localizes on $5$ sites:
$$\liminf_{n\rightarrow\infty } \frac{L_n\wedge R_n}{n}=0 \qquad \hbox{ but }\qquad  \limsup_{n\rightarrow\infty } \frac{L_n\wedge R_n}{n}>0.$$

Finally, let us mention that the functions $\varepsilon(n)$ and
$\varepsilon'(n)$ can also be explicitly computed for particular
weights sequences. For example, for $w(n)=n\log \log n$, using
\eqref{eqW} and similar arguments as in the proof of \eqref{Tx} , we
find that
$$W(\varepsilon(n))\equiv Y_n^+(x+1) \equiv \sum_{k=0}^{n/2}\frac{p_k}{w(k)}
 \quad \hbox{ and }\quad W(\varepsilon'(n))\equiv  W(n/2)- W(\varepsilon(n))$$
where $p_k$  denotes the probability of the walk to jump to site
$x+1$ at its $k$-th visit to $x+2$ (\emph{i.e.} $p_k\sim L_{2k}/k$).
Thus, after some (rather tedious) calculations, we deduce that, on
the event $\delta:=\lim_{n\rightarrow \infty} W(R_n)-W(L_n) \in
(0,1)$, the asymptotic local time profile on $R'$ takes the form:
\bigskip
\begin{center}
\includegraphics[width=9.6cm,angle=0]{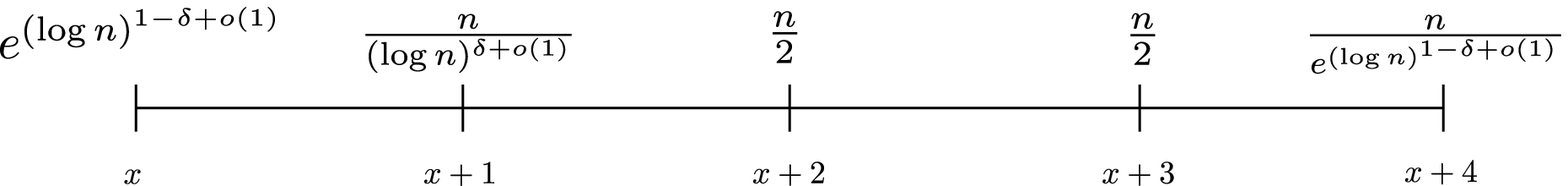}
\end{center}

\section{Appendix}\label{section5}

\begin{prop}
\label{alphac_prop} Let $w$ and $\widetilde{w}$ denote two
non-decreasing weight functions.
\begin{enumerate}
\item[(a)] For any $\lambda>0$, we have $\alpha_c(w) = \lambda \alpha_c(\lambda w)$
\textbf{(scaling)}.
\item[(b)] If $w \leq \widetilde w$, then $\alpha_c(w) \geq \alpha_c(\widetilde w)$
\textbf{(monotonicity)}.
\item[(c)] If $w\sim \widetilde w$, then $\alpha_c(w)=\alpha_c(\widetilde w)$ \textbf{(asymptotic
equivalence)}.
\end{enumerate}
\end{prop}
\begin{proof}
The scaling property (a) follows directly from the relation $\lambda
I_\alpha(\lambda w) = I_{\lambda \alpha}(w)$. We now prove (b). For
$x\ge 0$ and $\alpha>0$, set $u(x,\alpha) :=W^{-1}(W(x)+\alpha)$ and
define $\tilde{u}(x,\alpha)$ similarly for $\tilde{w}$. We have,
\begin{equation}\label{equ}
\int_x^{u(x,\alpha)} \frac 1 {w(t)}\, dt = \int_x^{\widetilde
u(x,\alpha)} \frac 1 {\widetilde w(t)}\, dt =\alpha.
\end{equation}
 When $w\le
\widetilde w$, the equality above implies that $u(x,\alpha)\le
\widetilde u(x,\alpha)$ for all $x$ and all $\alpha>0$. Since $w$ is
non-decreasing we get $w( u(x,\alpha)) \le \widetilde w(\widetilde
u(x,\alpha))$, hence $I_\alpha(w) \ge I_\alpha(\widetilde w)$. This
establishes (b).

 Suppose now that $w$ and
$\tilde{w}$ are two weight functions such that $w(x) = \tilde{w}(x)$
for all $x$ larger than some $x_0$. Then, in view of \eqref{equ}, we
see that $u(x,\alpha)=\tilde{u}(x,\alpha)$ for all $\alpha >0$ and
all $x\geq x_0$. Hence $\alpha_c(w) = \alpha_c(\tilde{w})$.  This
shows that $\alpha_c(w)$ does not depend upon the values taken by
$w$ on any compact interval $[0,x_0]$. Thus (c) follows directly
from (a) and (b).
\end{proof}

\begin{rem} Theorem \ref{4probapositive} states that when
$\alpha_c(w)$ is finite and non-zero, the walk localizes on either
$4$ or $5$ sites. It would certainly be interesting to estimate the
probability of each of these events. This seems a difficult
question. Let us remark that  these probabilities are not (directly)
related  to $\alpha_c(w)$. Indeed, for any $\lambda>1$, the weight
functions $w$ and $\lambda w$ define the same VRRW yet, $\alpha_c(w)
\neq \alpha_c(\lambda w)$.
\end{rem}

\begin{proof}[Proof of Proposition \ref{alphac_cor}]
We just need to check that $\alpha_c(x\log\log x) = 1$. For $x\geq
1$, set
\begin{equation*}
w(x): = x(1+L(\log x)),
\end{equation*}
where $L$ denotes the Lambert function defined as the solution of
$L(x)e^{L(x)} = x$. Then it follows from elementary calculation that
\begin{equation*}
\frac{w(W^{-1}(x))}{w(W^{-1}(x+\alpha))} =  \frac{x^x(1+\log
x)}{(x+\alpha)^{x+\alpha}(1+\log(x+\alpha))}
\underset{x\to\infty}{\sim} \frac{e^{-\alpha}}{x^\alpha}.
\end{equation*}
Therefore $\alpha_c(w) = 1$. Using now the well known equivalence
$L(x) \sim \log(x)$, we conclude using (c) of Proposition
\ref{alphac_prop} that $\alpha_c(x\log\log x) = 1$.
\end{proof}

\begin{lem}
\label{Halphaw} Assume that $w$ is non-decreasing. We have
\begin{equation*}
\liminf_{x\to\infty} \frac{w(x)}{x} \geq \frac{1}{\alpha_c(w)}.
\end{equation*}
In particular, when $\alpha_c(w) < \infty$ then $\sum 1/w(n)^2 <
\infty$ and if $\alpha_c(w) = 0$ then $w$ has super-linear growth.
\end{lem}
\begin{proof}
In view of the scaling property $\lambda I_\alpha(\lambda w) =
I_{\lambda \alpha}(w)$, we just need to prove that $\liminf w(x)/x <
1$ implies $I_{1}(w) = \infty$. Thus, let us assume that for some
$\varepsilon>0$, there exist arbitrarily large $x$ such that
$w(x)/x \leq 1-\varepsilon$. Then, for such an $x$ and for $y \leq
\varepsilon x$, we have, since $w$ is non-decreasing,
$$
W(y + (1-\varepsilon)x) - W(y) =
\int_{y}^{y+(1-\varepsilon)x}\frac{dz}{w(z)} \geq
\frac{(1-\varepsilon)x}{w(x)} \geq 1,
$$
which we can rewrite as
$$
\frac{1}{W^{-1}(W(y) + 1)} \geq \frac{1}{y+(1-\varepsilon)x}.
$$
Thus,
\begin{equation}\label{apen1}
\int_{\varepsilon x /2}^{\varepsilon x}\frac{1}{w(W^{-1}(W(y) + 1))}
\geq \frac{\varepsilon}{2-\varepsilon}.
\end{equation}
Since there exist arbitrarily large $x$ such that \eqref{apen1}
holds, we conclude that $I_{1}(w) = \infty$.

\end{proof}
\begin{rem}\label{remn!}
The previous lemma cannot be improved without additional assumptions
on $w$. Indeed, consider the weight function $w_0$ defined by
$w_0(x) = (n!)^{2}$ for $x \in [((n-1)!)^2,(n!)^2)$, $n\in \N^*$. It
is easily seen that
$$
\liminf_{x\to\infty} \frac{w_0(x)}{x} = 1 = \frac{1}{\alpha_c(w_0)}.
$$
This provides an example of a weight function which does not
uniformly grow faster than linearly and yet for which the VRRW
localizes on $4$ sites with positive probability. On the other hand,
if $w$ is assumed to be regularly varying, then using similar
arguments to those in the proof above, one can check that the
finiteness of $\alpha_c(w)$ implies $\lim_\infty w(x)/x = \infty$.
\end{rem}

\begin{lem}\label{lemaa} Assume that $w$ is non-decreasing and $\alpha_c(w) <
\infty$, for any $0 < \delta < \delta'$, we have
\begin{equation}\label{almostlast}
\liminf_{x\to\infty}\frac{W^{-1}(W(x) + \delta')}{W^{-1}(W(x) +
\delta)} \geq e^{\frac{\delta'-\delta}{\alpha_c(w)}}.
\end{equation}
Furthermore, for $\delta > \alpha_c(w)$,
\begin{equation}\label{andthelast}
\lim_{x\to\infty}\frac{x}{W^{-1}(W(x)+\delta)} = 0.
\end{equation}
As a consequence:
\begin{itemize}
\item[(a)] For any $\delta > \alpha_c(w)$ and any $c\in\R$,
$$
\sum_n^\infty \frac{1}{w(W^{-1}(W(n)+\delta)- c)} < \infty.
$$
\item[(b)] If $\alpha_c(w) = 0$, then, for any $\delta , \gamma > 0$,
$$
\sum_n^\infty \frac{1}{w(\gamma W^{-1}(W(n)+\delta))} < \infty.
$$
\end{itemize}
\end{lem}
\begin{proof}Recall the notation $u(x,\delta):= W^{-1}(W(x)+\delta)$. We
have
$$
\int_{u(x,\delta)}^{u(x,\delta')}\frac{ds}{w(s)} = \delta'-\delta.
$$
Using Lemma \ref{Halphaw} and the fact that $u(x,\delta)$ tends to
infinity as $x$ goes to infinity. we get, for any
$\alpha>\alpha_c(w)$,
$$
\liminf_{x\to
\infty}\log\left(\frac{u(x,\delta')}{u(x,\delta)}\right) =
\liminf_{x\to \infty} \int_{u(x,\delta)}^{u(x,\delta')}\frac{ds}{s}
\geq \lim_{x\to \infty}
\int_{u(x,\delta)}^{u(x,\delta')}\frac{ds}{\alpha w(s)} =
\frac{\delta'-\delta}{\alpha}
$$
which yields \eqref{almostlast}. Assertion (a) now follows from
\eqref{almostlast} noticing that, for $\alpha_c(w) < \gamma <
\delta$, we have $W^{-1}(W(n)+\delta)- c \geq W^{-1}(W(n)+\gamma)$
for all $n$ large enough. The proof of Assertion (b) is similar.

It remains to prove \eqref{andthelast}. Let $\delta > \alpha_c(w)$
and pick $\varepsilon >0$ small enough such that $\alpha_c(w) <
\delta - \varepsilon$. Assume by contradiction that, for some $A>0$,
we can find $x_0$ arbitrarily large such that $u(x_0,\delta) \leq A
x_0$.  Then, for all $x<x_0$,
\begin{equation*}
\varepsilon =
\int_{u(x,\delta-\varepsilon)}^{u(x,\delta)}\frac{dy}{w(y)} \leq
\int_0^{A x_0}\frac{dy}{w(u(x,\delta-\varepsilon))} = \frac{A
x_0}{w(u(x,\delta-\varepsilon))}
\end{equation*}
which, in turn, implies
\begin{equation*}
\int_{x_0/2}^{x_0}\frac{dx}{w(u(x,\delta-\varepsilon))} \geq
\frac{\varepsilon}{2A}
\end{equation*}
and contradicts the fact that $I_{\delta-\varepsilon}(w) < \infty$.
\end{proof}

\begin{lem}
\label{teclem} Assume that $w$ is non-decreasing. For any
$\beta<\alpha_c(w)$, we have
$$\sum_{n}^\infty \frac 1{w(n+W^{-1}(W(n)+\beta))} = \infty.$$
\end{lem}
\begin{proof}
Choose $\alpha\in (\beta,\alpha_c(w))$. Since $w$ is non-decreasing,
for any $t\ge 0$ and any $m\le n$, we have
\begin{equation}\label{lasteq}
u(n,t)-u(m,t)\ge n-m.
\end{equation}
Assume now that, for some large $n$, we have $n+u(n,\beta)\ge
u(n,\alpha)$. Then, necessarily, there exists $k\in
[u(n,\beta),u(n,\alpha)]$, such that $w(k)\le n/(\alpha-\beta)$. In
particular, since $w$ is non-decreasing, $w(u(n,\beta))\le
n/(\alpha-\beta)$. Moreover, using \eqref{lasteq}, we get
$m+u(m,\beta)\le u(n,\beta)$ for all $m\le n/2$. Thus we also have
$w(m+u(m,\beta))\le n/(\alpha-\beta)$, for all $m\le n/2$. It
follows that
$$\sum_{m=n/4}^{n/2} \frac 1{w(m+u(m,\beta))} \ge \frac{\alpha-\beta}{4}.$$
Therefore if $n+u(n,\beta)\ge u(n,\alpha)$, for infinitely many $n$,
the desired result follows. Conversely, if the inequality above
 holds only for finitely many $n$, then because $w$ is
non-decreasing, the result follows as well from the fact that
$I_\alpha(w)= \infty$.
\end{proof}

\vspace{0.2cm} \textit{Acknowledgments: This work started following
a talk given by Vlada Limic at the conference ``Marches
al\'eatoires, milieux al\'eatoires, renforcement'' in Roscoff
(2011) who shared some questions addressed to her by B\'alint T\'oth. We thank them both for having motivated this work. }

\bigskip

\footnotesize{This project is supported by ANR MEMEMO II, n°2010 BLAN 0125 02.}


\begin{thebibliography}{99}

\bibitem{BSS} A.-L. Basdevant, B. Schapira and A. Singh. Localization results for one-dimensional vertex reinforced random
walk with sub-linear reinforcement. \emph{In preparation}.

\bibitem{BT} M. Bena\"im and P. Tarr\`es. Dynamics of Vertex-Reinforced Random Walks. To appear in \emph{Ann. Probab}.

\bibitem{Ch} L.H.Y. Chen. A short note on the conditional Borel-Cantelli lemma. \emph{Ann. Probab.} 6(4):699--700, 1978.

\bibitem{CD} D. Coppersmith D and P. Diaconis. Random walk with reinforcement. Unpublished, 1987.

\bibitem{Dav} B. Davis. Reinforced random walk. \emph{Probab. Theory Related Fields} 84(2):203--229, 1990.

\bibitem{ETW1} A. Erschler, B. T\'oth and W. Werner. Some locally self-interacting walks on the integers. Preprint, arXiv:1011.1102.

\bibitem{ETW2} A. Erschler, B. T\'oth and W. Werner. Stuck Walks. To appear in \emph{Probab. Theory Related Fields}.

\bibitem{HH} P. Hall, C.C. Heyde. Martingale limit theory and its application. Probability and Mathematical Statistics. \emph{Academic Press, Inc.
[Harcourt Brace Jovanovich, Publishers], New York-London}, 1980.

\bibitem{P} R. Pemantle. Vertex-reinforced random walk. \emph{Probab. Theory Related Fields} 92(1):117--136, 1992 .

\bibitem{P2} R. Pemantle. A survey of random processes with reinforcement. \emph{Probab. Surv.} 4:1--79 (electronic), 2007.

\bibitem{PV}  R. Pemantle and S. Volkov. Vertex-reinforced random walk on $\Z$ has finite range. \emph{Ann. Probab.} 27(3):1368--1388, 1999.

\bibitem{Sch} B. Schapira. A $0-1$ law for Vertex Reinforced Random Walk on $\Z$ with weight of order $k^\alpha$, $\alpha<1/2$. Preprint, arXiv:1107.2505.

\bibitem{Sel} T. Sellke. Reinforced random walk on the $d$-dimensional integer lattice. \emph{Markov Process. Related Fields} 14(2):291--308, 2008.

\bibitem{T1} P. Tarr\`es. Vertex-reinforced random walk on $\Z$ eventually gets stuck on five points. \emph{Ann. Probab.} 32(3B):2650--2701, 2004.

\bibitem{T2} P. Tarr\`es. Localization of reinforced random walks. Preprint, arXiv:1103.5536.

\bibitem{V1} S. Volkov. Vertex reinforced random walk on arbitrary graphs, \emph{Ann. Probab.} 29(1):66--91, 2001.

\bibitem{V2} S. Volkov. Phase transition in vertex-reinforced random walks on $\Z$ with non-linear reinforcement. \emph{J. Theoret. Probab.} 19(3): 691--700, 2006.

\end{thebibliography}
\end{document}